\documentclass[12pt]{amsart}
\usepackage{amsfonts,yfonts,amsmath,amsthm,amssymb,bbm}
\usepackage[all]{xy}
\newcommand{\cts}[2]{\mathcal{A}^0(#1,#2)}
\newcommand{\tmop}[1]{\ensuremath{\operatorname{#1}}}
\newcommand{\assign}{:=}
\newcommand{\rl}{{\mathbbm{R}}}
\newcommand{\cx}{{\mathbbm{C}}}
\newcommand{\id}{{\mathbbm{I}}}
\newcommand{\m}{{\mathcal{M}}}
\newcommand{\dbar}{\overline{\partial}}
\newcommand{\Db}[1]{\frac{\partial{#1}}{\partial\overline{z}}}
\newcommand{\Dz}[1]{\frac{\partial{#1}}{\partial z}}
\newcommand{\Dx}[1]{\frac{\partial{#1}}{\partial x}}
\newcommand{\Dy}[1]{\frac{\partial{#1}}{\partial y}}
\newcommand{\abs}[1]{\left|{#1}\right|}
\newcommand{\norm}[1]{\left\|{#1}\right\|}
\newcommand{\hol}[2]{\mathcal{H}\left({#1},{#2}\right)}
\newcommand{\psh}{plurisubharmonic~}

\newtheorem{theorem}{Theorem}
\newtheorem{lemma}{Lemma}[section]
\newtheorem{prop}[lemma]{Proposition}

\newtheorem{definition}[lemma]{Definition}

\newtheorem{obs}[lemma]{Observation}
\title[Arcs and approximation]{Coordinate neighborhoods of arcs
 and the approximation of maps into (almost) complex manifolds}

\author{Debraj Chakrabarti}

\address{Department of Mathematics\\ University of Wisconsin-Madison}
 \email{dchakrab@math.wisc.edu}


\begin{document}
\begin{abstract}
We study the approximation of $J$-holomorphic maps continuous to the
boundary  from a domain in $\cx$ into an almost complex manifold by
maps $J$-holomorphic to the boundary, giving partial results in the
non-integrable case. For the integrable case, we study arcs in
complex manifolds and establish the existence of neighborhoods
biholomorphic to open sets in Euclidean space for several classes of
arcs. As an application we obtain $\mathcal{C}^k$ approximation of
holomorphic maps continuous to the boundary into complex manifolds
by maps holomorphic to the boundary, provided the boundary is nice
enough.
\end{abstract}
\maketitle
\section{Introduction}

This paper is divided into three sections, which though mostly
independent of each other, are devoted to the study of the following
question: let $\Omega\Subset\cx$, and let $(X,J)$ be an almost
complex manifold. Suppose that $f$ is a continuous map from the
closure $\overline{\Omega}$ to $X$ which is $J$-holomorphic on
$\Omega$. Can we approximate  $f$ by maps $J$ holomorphic on
(shrinking) neighborhoods of $\overline{\Omega}$?

In Section \ref{ac}  we give some conditions under which such a map
$f$ can be approximated by $J$-holomorphic maps in a neighborhood of
$\overline{\Omega}$. Unfortunately this involves smoothness
assumptions on the boundary $\partial\Omega$ and on $f$ as well (see
Theorem~\ref{almostcomplex} below.)

As might be expected,  when the almost complex structure $J$ is
integrable, we can say much more, since we have the tools of Complex
Analysis at our disposal. In fact, recently Drinovec-Drnov\u{s}ek
and Forstneri\u{c} have proved  the following : Let $S$ and $X$ be
complex manifolds, and let $\Omega\Subset S$ be a strongly
pseudoconvex Stein domain with boundary of class $\mathcal{C}^k$
with  $k\geq 4$. Then every $\mathcal{C}^{k-2}$ map from
$\overline{\Omega}$ to $X$ can be approximated in the
$\mathcal{C}^{k-2}$ sense by maps which are holomorphic on
(shrinking neighborhoods of) $\overline{\Omega}$ (see Theorem 2.1 in
\cite{forst}). This is a consequence of the fact that the graph of
such maps have a basis of Stein neighborhoods (Theorem 2.6 in
\cite{forst}.) In a forthcoming paper, they have been able to drop
the smoothness assumption on the map to be approximated, and avoid
the construction of Stein neighborhoods by use of the theory of
sprays. In this paper, we consider the case in which the source
manifold $S$ is the complex plane $\cx$, for which we give a proof
along completely different lines, and in the way obtain some results
of interest on their own.

Section \ref{coordinate} is devoted to the study of arcs (injective
continuous maps from the interval) in complex manifolds. We ask the
following question : under what circumstances does such an arc have
a coordinate neighborhood, i.e. a neighborhood biholomorphic to an
open set in Euclidean space $\cx^n$?  For a real analytic arc
$\alpha$ embedded in a complex manifold $\m$ (i.e., for each
$t\in[0,1]$ we have $\alpha'(t)\not=0$), it is an old result of
Royden that a coordinate neighborhood exists (see \cite{royden}).
 We show that embedded $\mathcal{C}^2$ arcs as
well as $\mathcal{C}^1$ arcs with some additional conditions have
coordinate neighborhoods (see Proposition~\ref{smoothnbds} and
Proposition~\ref{c1}.)For our application it is important not to
restrict attention to smooth arcs alone. We consider a class of
non-smooth arcs with finitely many non-smooth points which we call
mildly singular (see Definition~\ref{mildsingularitiesdefn} below.)
We show that such arcs have coordinate neighborhoods
(Theorem~\ref{singularities}).

 As an easy consequence of the
results of Section \ref{coordinate}, in Section \ref{approx} we
obtain the following result, a special case of the results of
Drivonec-Drov\u{s}nek and Forstneri\u{c} (with slightly weaker 
hypothesis on than theirs on the boundary):\\
{\em Let $k\geq 0$ be an integer and let  $\m$ be an arbitrary
complex manifold. Let $f$ be a $\mathcal{C}^k$ map from
$\overline{\Omega}$ into $\m$, where the open set $\Omega\Subset\cx$
is  bounded by finitely many Jordan curves, which are further
assumed to be $\mathcal{C}^1$ if $k\geq 1$. If $f$ is holomorphic on
$\Omega$, it can be approximated in the $\mathcal{C}^k$ topology by
holomorphic maps from (neighborhoods of) $\overline{\Omega}$ into
$\m$.}

This paper is based on the author's Ph. D. thesis (\cite{thesis}).
He would like to take this opportunity to express his deepest
gratitude to his advisor, Prof. Jean-Pierre Rosay. Without his
constant encouragement and help  neither the thesis and  nor this
paper would ever have been written.

\section{Maps into almost complex manifolds}

\label{ac} We begin by introducing some notation. For a compact
$K\subset\rl^N$, an integer $k\geq 0$ and $0<\theta<1$, let
$\mathcal{C}^{k,\theta}(K)$ denote the Lipschitz space of order
$k+\theta$, denoted by $\tmop{Lip}(k+\theta,K)$ in the text
\cite{stein}, where it is defined as a Banach space  of $k$-jets,
with $k$-th derivatives H\"{o}lder continuous with exponent
$\theta$.  If the set $K$ is nice, for example the closure of a
smooth domain (which will always be the case in our applications), we
can identify $\mathcal{C}^{k,\theta}(K)$ with the space of those $k$
times differentiable functions on $K$ all whose partial derivatives
of order $k$ are  H\"{o}lder continuous with exponent $\theta$. The
following remarkable fact is proved in \cite{stein}(p.177, Theorem
4.)
\begin{lemma}\label{genunifextension}
    Given $N,k\in \mathbbm{N}$ and $0<\theta<1$, there is a constant $C$ with the following property.
    Given any compact $K\subset \rl^N$, there is a linear extension operator
    $E:\mathcal{C}^{k,\theta}(K)\rightarrow \mathcal{C}^{k,\theta}(\rl^N)$ such that
    $\norm{E}_{\tmop{op}} < C$.
\end{lemma}

For a Riemannian manifold $(X,g)$, we define the space
$\mathcal{C}^{k,\theta}(K,X)$ of Lipschitz maps in the obvious way
using local charts on $X$, and this space has a natural structure of a
metric space. Observe however that the topology on
$\mathcal{C}^{k,\theta}(K,X)$ is not dependent on the choice of the
metric $g$.

Let $(X,J)$ be an almost complex manifold, where we assume that the
almost complex structure $J$ is of class $\mathcal{C}^{k,\theta}$
for some $k\geq 1$ and $0<\theta<1$.  For compact $K\subset\cx$, we
denote by $\mathcal{H}_J(K,X)$ the space of $J$-holomorphic maps from
$K$ to $X$. The map  $f:K\rightarrow
X$ belongs to $\mathcal{H}_J(K,X)$, iff for some open $U_f\supset K$, 
the map $f$  extends $J$-holomorphically to $U_f$.  It is well known that
$\mathcal{H}_J(K,X)\subset \mathcal{C}^{k+1,\theta}(K,X)$.  Let
$\mathcal{A}^{k,\theta}_J(K,X)$ denote the closed subspace of
$\mathcal{C}^{k,\theta}_J(K,X)$ consisting of maps $f$ which are
$J$-holomorphic on the interior of  $K$.

We can now state the following:

\begin{theorem}\label{almostcomplex}
Let $\left(X,J\right)$ be an almost complex manifold with the
structure $J$ of class $\mathcal{C}^{k,\theta}$, where $k\geq 1$ and
$0<\theta<1$, and let $\Omega$ be an open set in $\cx$ with
$\mathcal{C}^1$ boundary. Then in the metric space
$\mathcal{C}^{k,\theta}(\overline{\Omega},X)$, the set
$\mathcal{H}_J(\overline{\Omega},X)$ is dense in the set
$\mathcal{A}^{k+1,\theta}_J(\overline{\Omega},X)$.
\end{theorem}

First we prove a slightly stronger
 version (Proposition~\ref{r2n})  for  $X=\rl^{2n}$, and then reduce
 the general case of an a.c. manifold $X$ to it.
\begin{lemma}\label{lhg}
    Let $k\geq 1$ be an integer,  $\omega\Subset\cx$  an open set, and  $B$ a $2n\times 2n$
        real matrix of $\mathcal{C}^{k-1,\theta}$
     functions on $\overline{\omega}$. Let $L$ denote the differential operator given by
\[      L(h) = \Db{ h} + B h, \]
    mapping $\mathcal{C}^{k,\theta}(\overline{\omega})$ into $\mathcal{C}^{k-1,\theta}(\overline{\omega})$.
    There exists a constant $C_0$ such that for any open subset $W\subset\omega$, and any
    $g\in\mathcal{C}^{k-1,\theta}(\overline{W})$ there exists $h\in\mathcal{C}^{k,\theta}(\overline{W})$ such
    that, on $\overline{W}$,
    \begin{equation}\label{lhgeq}
            Lh =g
    \end{equation}
    and
\[      \norm{h}_{\mathcal{C}^{k,\theta}(\overline{W})}
            \leq C_0 \norm{g}_{\mathcal{C}^{k-1,\theta}(\overline{W})}.\]
\end{lemma}
\begin{proof}
Fix $R>0$ such that $\omega\subset\Delta_R=
\left\{z\in\cx\colon\abs{z}<R\right\}$. By Theorem A2 of
\cite{rosayivashkovich}, one can solve the equation
\begin{equation}\label{ivashkovichequation}
  \Db{\tilde{h}} + B \tilde{h} = \tilde{g}
  \end{equation} on $\Delta_R$ with
$\norm{\tilde{h}}_{\mathcal{C}^{k,\theta}(\Delta_R)}\leq K_R
\norm{\tilde{g}}_{\mathcal{C}^{k-1,\theta}(\Delta_R)}$, where $K_R$
is a constant depending only on $R$. (The proof in
\cite{rosayivashkovich} assumes that $k=1$, but it generalizes
immediately)

Let $C$ be the absolute constant provided by
Lemma~\ref{genunifextension} as an upper bound to the norm of linear
extension operators mapping $\mathcal{C}^{k-1,\theta}$ of a compact
subset of $\rl^2$ to $\mathcal{C}^{k-1,\theta}(\rl^2)$. Extend the
data $g$ and the coefficients $B$ of equation~(\ref{lhg}) from
$\overline{W}$ to $\tilde{g}$ and $\tilde{B}$ defined on $\rl^2$,
with $\norm{g}_{\mathcal{C}^{k-1,\theta}(\rl^2)}\leq C
\norm{g}_{\mathcal{C}^{k-1,\theta}(\overline{W})}$, and similarly
for $B$. We now solve equation \ref{ivashkovichequation} with
estimates as mentioned above, and  set $h$ to be the restriction of
$\tilde{h}$ to $\overline{W}$.
\end{proof}

We now prove a version of Theorem~\ref{almostcomplex} for
$X=\rl^{2n}$.
\begin{prop}
\label{r2n}
    Let  $\Omega\Subset\cx$ be an open set and let $U$ be an open
    neighborhood of $\overline{\Omega}$.
     Let $J$ be an almost complex structure
    of class $\mathcal{C}^{k,\theta}$ on  $\rl^{2n}$, with $k\geq
    1$.  Let $\beta$ be such that $\theta<\beta<1$.
    Suppose that   the $\mathcal{C}^{k  ,\beta}$ map
    $f:\overline{U}\rightarrow \rl^{2n}$ is such that
    \begin{itemize}
        \item $f|_\Omega$ is $J$-holomorphic, and
        \item $J|_{f(\Omega)}= J_{\tmop{st}}$, the standard complex structure of $\rl^{2n}=\cx^n$.
    \end{itemize}
    Then $f$ can be approximated uniformly on $\overline{\Omega}$ by $J$-holomorphic maps.
\end{prop}
\begin{proof}
Suppose we are given $\epsilon_0>0$. We want to find a neighborhood
$\Omega_{\epsilon_0}$ of $\overline{\Omega}$, and a $J$-holomorphic
$u$ from $\overline{\Omega}_{\epsilon_0}$ into $R^{2n}$  such that
$\norm{u-f}_{\mathcal{C}^{k,\theta}(\overline{\Omega}_{\epsilon_0})}<\epsilon_0$.

We set
\[ \Db{u} =\frac{1}{2}\left(\Dx{u} + J_{\tmop{st}}\Dy{u}\right) \]
and
\[ \Dz{u} =\frac{1}{2}\left(\Dx{u} - J_{\tmop{st}}\Dy{u}\right) \]
since $J_{\tmop{st}}$ corresponds to multiplication by $i$ in the
identification of $\rl^{2n}$ with $\cx^n$. We know that, provided
that  $J+J_{\tmop{st}}$ is invertible, the
 condition that a map $u$ from
some subset of $\cx$ into $(\rl^{2n},J)$ is $J$-holomorphic is given
by
\[ \Phi ( u ) \assign \Db{u} + Q ( u )\Dz{u}=0\]
where $Q ( u )$ is a $2n \times 2n$  matrix, given by $Q ( u ) = [ J
( u ) + J_{\text{st}} ]^{- 1} [ J ( u ) - J_{\tmop{st}} ]$.  Since
$J=J_{st}$ on the range of $f$, for maps $u$ sufficiently close to
$f$, we have $J(u)\approx J_{\tmop{st}}$, so this equation
determines the $J$-holomorphy of $u$ for such maps.

We will think of $\Phi$ to be a map from $\mathcal{C}^{k, \theta} (
\overline{U} )$ to $\mathcal{C}^{k-1,\theta} (\overline{U})$ . Its
derivative is given by:
\begin{eqnarray}
\label{derivative} \Phi' ( u ) h   & =&    \frac{\partial
h}{\partial \bar{z}} + Q' ( u ) h
            {\frac{\partial u}{\partial z}} + Q ( u ) {\frac{\partial
            h}{\partial z}}\nonumber \\
        &= &    \left\{ \frac{\partial h}{\partial \bar{z}} + A ( u ) h
            \right\} + Q ( u ) {\frac{\partial h}{\partial z}}\nonumber\\
        &=&     L_u h + R_u h.
\end{eqnarray}
Observe that $A$ and $Q$ are $2n\times 2n$ matrices with entries in
$\mathcal{C}^{k-1,\theta}(\overline{U})$ and
$\mathcal{C}^{k,\theta}(\overline{U})$ respectively. Since we can
easily show that the assignments $u\mapsto (h\mapsto A(u)h)$ and
$u\mapsto (h\mapsto Q(u)h_z)$ are continuous from
$\mathcal{C}^{k,\theta}(\overline{U})$ into the Banach space of
operators
 $BL\left(\mathcal{C}^{k,\theta}(\overline{U},\rl^{2n}),
\mathcal{C}^{k-1,\theta}(\overline{U},\rl^{2n})\right)$,  it follows
that $\Phi$ is $\mathcal{C}^1$.

We claim the following~: {\em There is an open
$W\supset\overline{\Omega}$ such that $\Phi'(f)$ is surjective
from $\mathcal{C}^{k,\theta}(\overline{W})$ to $\mathcal{C}^{k-1,\theta}(\overline{W})$.}\\
To prove the claim,  observe that $Q(f)\in \mathcal{C}^{k,\theta}$,
therefore, a fortiori $Q(f)$ is in $\mathcal{C}^k$. We can choose
$W\supset\Omega$ so small that
$\norm{Q(f)}_{\mathcal{C}^k(\overline{W})}$ is small. Since the
boundary $\partial\Omega$ of the set $\Omega$ is $\mathcal{C}^1$ by
hypothesis, we can also choose the $W$ such that the
$\mathcal{C}^{k-1,\theta}$ norm is dominated by the $\mathcal{C}^k$
norm. Therefore, we can find a $W$ such that
\[ \norm{Q(f)}_{\mathcal{C}^{k-1,\theta}(\overline{W})} <\frac{1}{2C_0 K}, \]
where $C_0$ is the constant in Lemma~\ref{lhg} above. Therefore we
have for $h$ in $\mathcal{C}^{k,\theta}(\overline{W})$:
\begin{eqnarray*}
\norm{R_fh}_{\mathcal{C}^{k-1,\theta}(\overline{W})}
    &=& \norm{Q(f)\frac{\partial h}{\partial\overline{z}}}_{\mathcal{C}^{k-1,\theta}(\overline{W})}\\
    &\leq &\norm{Q(f)}_{\mathcal{C}^{k-1,\theta}(\overline{W})}
        \norm{\frac{\partial h}{\partial\overline{z}}}_{\mathcal{C}^{k-1,\theta}(\overline{W})}\\
    &\leq & \norm{Q(f)}_{\mathcal{C}^{k-1,\theta}(\overline{W})}
        \norm{h}_{\mathcal{C}^{k,\theta}(\overline{W})}\\
    & < & \frac{1}{2C_0 K}\norm{h}_{\mathcal{C}^{k,\theta}(\overline{W})},
\end{eqnarray*}
so that $\norm{R_h}_{\tmop{op}}< \frac{1}{2C_0 K}$. Therefore,
$\Phi'(f)$ is a small perturbation of a surjective linear map, and
standard methods based on iteration shows that it is  surjective as
a map from $\mathcal{C}^{k,\theta}(\overline{W})$ to
$\mathcal{C}^{\theta}(\overline{W})$, and the equation $\Phi'(f)h=g$
can be solved with
$\norm{h}_{\mathcal{C}^{k,\theta}(\overline{W})}\leq 2C_0 K
\norm{g}_{\mathcal{C}^{\theta}(\overline{W})}.$

 Since $\Phi'(f)$ is surjective, and $\Phi$ is $\mathcal{C}^1$, we see that
there is a small ball around $f$ which is mapped surjectively by
$\Phi$ onto a ball around $\Phi(f)$. Therefore, given $\epsilon>0$
there is a $\delta>0$ such that   such that if $g$ in
 $\mathcal{C}^{k-1,\theta}(\overline{W})$ is such that
$\norm{g}_{\mathcal{C}^{k-1,\theta}(\overline{W})}<\delta$, then we
can solve the equation
\[ \Phi(f+r) =\Phi(f)+ g \]
for an $r\in \mathcal{C}^{k,\theta}(\overline{W})$ such that
 $\norm{r}_{\mathcal{C}^{k,\theta}(\overline{W})}<\epsilon$.

Now we fix $\epsilon=\epsilon_0$ ( where $\epsilon_0$ is as in the
beginning of this proof) and denote by $\delta_0$ the corresponding
$\delta$. Let $C$ be a uniform bound for linear extension operators
from $\mathcal{C}^{k-1,\theta}(\overline{V})$ to
$\mathcal{C}^{k-1,\theta}(\cx)$ for any open subset $V$ of $\cx$
(see Lemma~\ref{genunifextension}) and let
$\delta_1=\frac{\delta_0}{C}$. Since by hypothesis
$f\in\mathcal{C}^{1,\beta}$, it follows that
$\Phi(f)\in\mathcal{C}^\beta$. We now use the hypothesis that
$\beta>\theta$. Since $\Phi$ vanishes on $\Omega$, in a small enough
neighborhood of $\overline{\Omega}$, we have that $\norm{\Phi(f)}$
is small in the $\mathcal{C}^{k-1,\theta}$ sense. Let
$\Omega_{\epsilon_0}$ be a neighborhood of $\overline{\Omega}$ such
that we have
$\norm{\Phi(f)}_{\mathcal{C}^{k-1,\theta}(\Omega_{\epsilon_0})}<\delta_1$.
Denote by $g$ the map $-\Phi(f)|_{\Omega_{\epsilon_0}}$. Using a
linear extension operator, we extend $g$ to a function $\tilde{g}$
in $\mathcal{C}^{k-1,\theta}(\overline{W})$, such that
$\norm{\tilde{g}}_{\mathcal{C}^{k-1,\theta}(\overline{W})}\leq
C\delta_1= \delta_0$. Therefore, the equation
$\Phi(f+r)=\Phi(f)+\tilde{g}$ can be solved with for an $r$ such
that $\norm{r}_{\mathcal{C}^{k,\theta}(\overline{W})} <\epsilon_0$.
Now, if we set $u=f+r$, we have that on $\Omega_{\epsilon_0}$ we
have $\Phi(u)=-g+g=0$, i.e. $u$ is $J$-holomorphic. Of course, we
have
\begin{eqnarray*}
\norm{u-f}_{\mathcal{C}^{k,\theta}(\overline{\Omega}_{\epsilon_0})}
&\leq&\norm{u-f}_{\mathcal{C}^{k,\theta}(\overline{W})}\\
&\leq&\norm{r}_{\mathcal{C}^{k,\theta}(\overline{W})}\\
&<&\epsilon_0.
\end{eqnarray*} \nopagebreak
\end{proof}
\subsection{The general case}

Now let $(X,J)$ be an almost complex manifold, with $J$ of class
$\mathcal{C}^{k,\theta}$, with $k\geq 1$. Let
$f\in\mathcal{A}^{k+1,\theta}_J(\overline{\Omega},X)$. To prove
Theorem~\ref{almostcomplex} we need to approximate $f$ in the
$\mathcal{C}^{k,\theta}$ topology on $\overline{\Omega}$ by
$J$-holomorphic maps.

We begin by making two observations. First, that it is no loss of
generality to assume that $f$ is  an embedding. This is because we
can replace $X$ by $\mathbbm{C} \times X$ and replace  $f$ by the
map $F : z \mapsto ( z, f ( z ) )$,  and obtain an approximation to
$F$, which we can subsequently project onto $X$. We will therefore
assume to begin with that $f$ is an embedding.

The next observation is that we can extend $f$ as a
$\mathcal{C}^{k+1,\theta}$ map to all of $\mathbbm{C}.$ Therefore we
will assume that $f$ is defined and is an embedding on some large
set $\overline{U}$ containing the set $\Omega$ compactly, and $f$ is
$J -$holomorphic on $\Omega$.

Let $n$ denote the complex dimension of the a.c. manifold $X.$ It is
easy to find $( n - 1 )$ smooth vector fields $Y_2, Y_3, \ldots$,
$Y_n$ on the embedded disc $f (\overline{U})$ such that for any
point $z$  the $\mathbbm{C}-$span of the vectors $\frac{\partial
f}{\partial x} ( z ), Y_2 ( f ( z ) ), \ldots, Y_n ( f ( z ) )$ in
the space $T_{f ( z )} X$ with respect to the complex structure
induced by $J ( f ( z ) )$ is the whole of $T_{f ( z )} X$. Now
consider the map from $\overline{U} \times \mathbbm{C}^{n - 1}
\subset \mathbbm{C}^n$ into $X$ given by
\begin{equation}
  \label{eq1} ( z_1, \ldots, z_n ) \mapsto \exp_{\sum^n_{j = 2} z_j Y_j ( f (
  z_1 ) )} ( f ( z_1 ) ) .
\end{equation}

There is a neighborhood of $\overline{U} \times \mathbbm{C}^{n - 1}$
in $\cx^n$ which is mapped by the map in equation (\ref{eq1})
diffeomorphically onto a neighborhood $\mathcal{U}$ of $f
(\overline{U})$ in $X$. Therefore the inverse $\varphi$ of the map
in equation (\ref{eq1}) is a system of coordinates on $\mathcal{U}$.
We note some properties of this coordinate map:
\begin{itemize}
    \item $\varphi$ is of class $\mathcal{C}^{k+ 1,\theta}$. Consequently, the smoothness of  $J$
    is preserved, i.e., the induced structure $J^\sharp$ on  $\rl^{2n}$ is still $\mathcal{C}^{k,\theta}$.
    \item The map $f$ is represented in these coordinates by
     \begin{equation}
            \label{eq2} \zeta \mapsto ( \zeta, \underset{n - 1}{\underbrace{0, \ldots,
            0}} ) \in \mathbbm{C}^n .
    \end{equation}
  \item   On the set $ \overline{\Omega} \times \{
  \mathbf{0^{n - 1} \}}$ we have that $J^\sharp = J_{\text{st}}$,  the standard almost
  complex structure of $\mathbbm{C}^n$.
\end{itemize}
The problem is therefore reduced to that considered in the first
section, and the approximation asserted in the theorem can be done.
\section{ Arcs in  Complex Manifolds}

\label{coordinate}
 We will denote by $\m$ a complex manifold of
complex dimension $n$ on which we impose a Riemannian metric $g$.
The actual choice of the metric does not affect any of our results.

An {\em arc} is an injective continuous map from the interval
$[0,1]$. We say that a $\mathcal{C}^1$ arc  $\alpha$ is {\em
embedded} if $\alpha'(t)\not=0$ for each $t$. For convenience of
exposition we introduce the following terminology:
\begin{definition}{\rm
\label{goodprojection} Let  $\alpha$ be an arc in $\m$. Let $\phi$
be a holomorphic submersion from a neighborhood of  $\alpha([0,1])$
in $\m$ into $\cx$. We will say that $\phi$ is a {\em good
submersion} for the arc $\alpha$ if $\phi\circ\alpha$ is a
$\mathcal{C}^1$ embedded arc in $\cx$}
\end{definition}
Clearly, a smooth (at least $\mathcal{C}^1$) arc which admits a good
submersion is embedded. Observe however that the definition does not
require the arc to be smooth. Indeed, the existence of good
submersions will serve as a convenient substitute for being embedded
when we consider non-smooth arcs.

First we generalize Royden's result on the existence of co-ordinate
neighborhoods of real analytic arcs to smooth arcs. The proof of
this result is based on a quantitative approximation of 
$\mathcal{C}^k$ arcs  by real analytic arcs (see
Lemma~\ref{analyticapprox} below.)

 Our first result is as follows:\\
 \begin{prop}
 \label{smootharc}
{ Let $k\geq 2$ and let $\alpha$ be an embedded $\mathcal{C}^k$ arc
in $\m$. Then there is a family $\{\phi_j\}_{j=1}^n$ of  $n$ good
submersions associated with $\alpha$, such that they form a
coordinate system in a neighborhood of the image of $\alpha$. }
\end{prop}

In particular, smooth arcs of class at least $\mathcal{C}^2$ have
coordinate neighborhoods. Also, a $\mathcal{C}^2$ arc is embedded
iff it has a good submersion.

We next consider  $\mathcal{C}^1$ arcs in $\m$. Unfortunately, in
this case, the approximation lemma \ref{analyticapprox} is not
strong enough to prove the existence of coordinate neighborhoods, if
we simply assume that $\alpha$ is embedded. However, we can prove
the
following:\\
\begin{prop}
\label{c1} Let $\alpha$ be a $\mathcal{C}^1$ arc in $\m$ which
admits a good submersion $\phi$. There is a coordinate neighborhood
$W$ of $\alpha([0,1])$ in $\m$  and a biholomorphic map
$(\phi_1,\ldots,\phi_n)$ from $W$ onto on open subset of $\cx^n$
such that $\phi_n=\phi|_W$.
\end{prop}

In other words, given a good submersion in a neighborhood of the
image of an arc in $\m$, we can find $n-1$ other functions such that
the $n$ functions together form a system of coordinates in a
neighborhood of $\alpha$. Observe that, for $j=1,\ldots,n-1$, after
replacing the function $\phi_j$ by the function $\phi_j+ K\phi$, for
large enough $K$, we may assume that each of the coordinate
functions $\phi_j$ is a good submersion, thus strengthening the
conclusion.

We now turn to non-smooth arcs. In view of the intended application
in the next section we introduce the following definition:
\begin{definition}\label{mildsingularitiesdefn}{\rm
Let $k\geq 1$. Suppose
$\alpha:[0,1]\rightarrow\m$ is an arc  such that\\
\noindent$\bullet$ $\alpha$ is  $\mathcal{C}^k$ outside a finite
subset  $P\subset[0,1]$, and\\
\noindent$\bullet$ $\alpha$ admits a good submersion $\phi$.

We will refer to such arcs $\alpha$ as {\em $\mathcal{C}^k$ arcs
with mild singularities} or {\em mildly singular arcs}}
\end{definition}

Our result concerning mildly singular arcs is as follows:
\begin{theorem}\label{singularities}
Let $\alpha$ be a $\mathcal{C}^3$ arc in $\m$ with mild
singularities, and let  $\phi$ be the associated good submersion.
Then the image $\alpha([0,1])$ has a coordinate neighborhood $W$ and
a coordinate map $(\phi_1,\cdots,\phi_n):W\rightarrow\cx^n$, with
$\phi_n=\phi|_W$.
\end{theorem}

\subsection{Approximation of smooth functions by real-analytic functions}

The following approximation lemma is required in the proof of
the fact that smooth arcs have coordinate neighborhoods.

\begin{lemma}
\label{analyticapprox} Let $\Gamma$ denote either the image of the
unit interval $[0,1]$ or of the unit circle $S^1$ under a
$\mathcal{C}^k$ embedding into $\cx$, where $k\geq 1$. Let $f$ be a
$\mathcal{C}^k$ function defined on $\Gamma$, and $\theta$ be such
that $0<\theta<1$. Then there is a constant $C>0$ and a
$\mathcal{C}^k$ extension of  $f$ to a neighborhood of $\Gamma$ such
that for small enough $\delta>0$ there is a holomorphic map
$f_\delta$ defined in the closed $\delta$ neighborhood
$\overline{B_\cx(\Gamma,\delta)}$ of $\Gamma$ such that
\begin{itemize}
\item if $\alpha$ and $\beta$ be nonnegative integers such that
$\alpha+\beta<k$, then for  $z\in \overline{B_\cx(\Gamma,\delta)}$we
have
\[ \abs{\left(\frac{\partial^{\alpha+\beta}}
{\partial z^\alpha \partial\overline{z}^\beta}\right)
\left(f(z)-f_\delta(z)\right)} < C
\delta^{k-\frac{1}{2}-(\alpha+\beta)}. \] Further,

\item $f_\delta$ is bounded independently of $\delta$ in the
$\mathcal{C}^{k-1,\theta}$ norm, or more precisely, we have
$\norm{f_\delta}_{\mathcal{C}^{k-1,\theta}
\left(\overline{B_\cx(\Gamma,\delta)}\right)} <C.$
\end{itemize}
\end{lemma}
\begin{proof}

We will denote by  $C$ any constant which is independent of
$\delta$.

For small $\delta>0$, let $\chi_\delta$ be a $\mathcal{C}^\infty_c$
cutoff on $\cx$ such that $\chi_\delta\equiv 1$ in a neighborhood of
$\overline{B_\cx(\Gamma,\delta)}$ and vanishes off the
$2\delta$-neighborhood ${B_\cx(\Gamma,2\delta)}$. We may choose the
$\chi_\delta$ that there is a constant $C$ (which of course depends
on $k$) such that for small $\delta$  and every pair of non-negative
integers $\alpha$ and $\beta$ such that $\alpha+\beta<k$ we have:
\[\abs{\left(\frac{\partial^{\alpha+\beta}}
{\partial z^\alpha \partial\overline{z}^\beta}\right) \chi_\delta(z)}
< \frac{C}{\delta^{\alpha+\beta}}.\]

Since $\Gamma$ is totally real,  using the Whitney Extension
Theorem, we can extend $f$ as a $\mathcal{C}^k$ function on $\cx$
such that the $\dbar$-derivative
 $\tfrac{\partial f}{\partial \overline{z}}$ vanishes to order $k-1$
 on $\Gamma$. Continuing to denote the extended function by $f$ and denoting by $\eta(z)$ the
distance from $z\in\cx$ to the set $\Gamma$, we see that
\[\abs{\left(\frac{\partial^{\alpha+\beta}}
{\partial z^\alpha \partial\overline{z}^\beta}\right)
 \left(\Db{f}\right)(z)}< {C}{\eta(z)^{k-1-(\alpha+\beta)}}.\]

Now we define (suppressing the dependence on $\delta$ in the
notation):
\[ \lambda_{(\alpha,\beta)}(z) =\frac{\partial^{\alpha+\beta}}
{\partial z^\alpha \partial\overline{z}^\beta}
\left(\chi_\delta\cdot\Db{f}\right)(z).\] Observe that
$\lambda_{(\alpha,\beta)}$ is supported in $B_\cx(\Gamma,2\delta)$
for every $\delta$ and we have $\abs{\lambda_{(\alpha,\beta)}}
=O\left(\delta^{k-1-(\alpha+\beta)}\right)$.

Let the function $u_\delta$ on $\cx$ be defined by:
\[
 u_\delta(z)\assign \frac{-1}{\pi z}*\lambda_{(0,0)}(z)
=\frac{-1}{\pi z}*\left(\chi_\delta(z)\cdot\Db{f}(z)\right).
\]

Then $f_\delta= f+u_\delta$ is clearly holomorphic on
$\overline{B_\cx(\Gamma,\delta)}$,
 and we have
\begin{eqnarray*}
\frac{\partial^{\alpha+\beta}}{\partial z^\alpha
\partial\overline{z}^\beta}
u_\delta(z)&=&\frac{-1}{\pi z}* \frac{\partial^{\alpha+\beta}}
{\partial z^\alpha
\partial\overline{z}^\beta}\left(\lambda_{(0,0)}(z)\right)
=\frac{-1}{\pi z}* \lambda_{(\alpha,\beta)}(z) \\
&=&-\frac{1}{\pi} {\int\int}_{B_\cx(\Gamma,2\delta)}
\frac{1}{z-\zeta}\cdot
\lambda_{(\alpha,\beta)}(\zeta)d\xi d\eta \hspace{1.5cm} (\zeta=\xi+\imath\eta)\\
&=& -\frac{1}{\pi}\left({\int\int}_{B_\cx(\Gamma,2\delta)\cap\{\zeta\colon\abs{\zeta-z}<
\sqrt\delta\}}+{\int\int}_{B_\cx(\Gamma,2\delta)\cap{\zeta\colon\{\abs{\zeta-z}\geq\sqrt\delta\}}}
 \right)\\
&=& -\frac{1}{\pi}\left(I_1 + I_2\right).
\end{eqnarray*}
We can now estimate:
\begin{eqnarray*}
\abs{I_1}&\leq& \frac{1}{\pi}
\norm{\lambda_{(\alpha,\beta)}}_{L^\infty}
{\int\int}_{\{\zeta\colon\abs{\zeta-z}<\sqrt\delta\}}
\frac{1}{\abs{z-\zeta}}d\xi d\eta\\
&\leq&\frac{1}{\pi}\cdot C\delta^{k-1-(\alpha+\beta)}\cdot
2\pi\sqrt\delta\\
&\leq &C\delta^{k-\frac{1}{2}- (\alpha+\beta)}
\end{eqnarray*}
and
\begin{eqnarray*}
\abs{I_2}&\leq& \frac{1}{\pi}
\norm{\lambda_{(\alpha,\beta)}}_{L^\infty}
{\int\int}_{B_\cx(\Gamma,2\delta)\cap\{\zeta\colon\abs{\zeta-z}\geq\sqrt\delta\}}
\frac{1}{\abs{z-\zeta}}d\xi d\eta\\
&\leq&\frac{1}{\pi}\cdot C\delta^{k-1-(\alpha+\beta)}\cdot
\frac{1}{\sqrt\delta}\cdot \mbox{Area}\left(B_\cx(\Gamma,2\delta)\right)\\
&\leq& C \delta^{k-1-(\alpha+\beta)}\cdot\frac{1}{\sqrt\delta} C \delta\\
&\leq & C\delta^{k-\frac{1}{2}- (\alpha+\beta)}
\end{eqnarray*}
which proves the first conclusion of the lemma. Moreover, it
immediately follows that for $f_\delta= f+u_\delta$ we have that
$\norm{f_\delta}_{\mathcal{C}^{k-1}} <C$ for some $C$ independent of
$\delta$. To complete the proof, it is sufficient to recall the well
known fact that for functions $v$ on $\cx$ supported in a fixed
compact set $E$, the assignment $v\mapsto \frac{1}{z}* v$ is
continuous from $\mathcal{C}(E)$ to $\mathcal{C}^{0,\theta}(E)$.
\end{proof}

\subsection{$\mathcal{C}^k$ embedded arcs, $k\geq 2$}
This section is devoted to the proof of Proposition~\ref{smootharc}.
We will need the following two lemmas:
\begin{lemma}\label{babysurjective}
For convenience, let $B = B_{\rl^N}(0,1)$, the $N$-dimensional unit
ball. Let $\Phi:\overline{B}\rightarrow\rl^N$ be a $\mathcal{C}^1$
map such that for some constant $C>0$,
\begin{itemize}
\item for each tangent vector $\mathbf{v}$, we have
$\norm{\Phi'(0){\mathbf v}}\geq C \norm{\mathbf v}$, and
\item for each $x\in \overline{B}$, we have
$\norm{\Phi'(x)-\Phi'(0)}_{\tmop{op}}< \frac{C}{2}.$
\end{itemize}
Then, $ \Phi(\overline{B})\supset
B_{\rl^N}\left(\Phi(0),\frac{C}{2}\right). $
\end{lemma}
\begin{proof}
After a translation and dilation, we can assume that $\Phi(0)=0$ and
$C=2$. Fix $x\in \overline{B}$ and let $u(t)=\Phi(tx)$.  We have:

\begin{eqnarray*}
\norm{\Phi(x)}&=& \norm{ \int_0^1u'(t) dt}
=\norm{ \int_0^1 \Phi'(tx)x dt}\\
&=& \norm{\int_0^1\Phi'(0)xdt +
 \int_0^1\left(\Phi'(tx)-\Phi'(0) \right)dt}\\
&\geq& \norm{\Phi'(0)x} -
\norm{\int_0^1\left(\Phi'(tx)-\Phi'(0) \right)xdt}\\
&\geq& 2\norm{x} -\norm{x} \geq  \norm{x}.
\end{eqnarray*}
 Let $\Sigma=B\cap \Phi(\overline{B})$. Then $\Sigma$ is  nonempty
($0\in\Sigma$) and closed in the relative topology of $B$. Now as
$\Phi$ is expanding, $B\cap \Phi(\overline{B})=B\cap \Phi({B})$.
Since $\Phi$ is a local diffeomorphism, it follows from this that
$\Sigma$ is open in $B$ as well, which implies by connectedness that
$\Sigma=B$, that is, $B\subset \Phi(\overline{B})$, which is the
required conclusion.
\end{proof}

\begin{lemma}\label{totallyreal}
    Let $\mathcal{S}$ be a sufficiently smooth compact totally real
    submanifold of  $\m$.
    Then, there is an $\eta>0$ such that the $\eta$ neighborhood
     $B_\m(\mathcal{S},\eta)$
    of  $\mathcal{S}$ is a Stein open subset of  $\m$
    and any continuous function on $\mathcal{S}$ can be uniformly approximated
    by the restrictions to $\mathcal{S}$ of functions holomorphic on this $\eta$
    neighborhood.

    ``Sufficiently smooth" in this context means $\mathcal{C}^s$, where $s$ is
    an integer which is at least 2  and greater than
    $\frac{1}{2}\mbox{dim}_\rl\mathcal{S}+1$.
\end{lemma}
In the case when $\mathcal{S}$ is  $\mathcal{C}^\infty$, this result
is due to Nirenberg and Wells (see \cite{nirenbergwells}, Theorem
6.1 and Corollary 6.2.) Since the submanifold $\mathcal{S}$ is of
class at least $\mathcal{C}^2$, the square of the distance to
$\mathcal{S}$ is strictly plurisubharmonic in a neighborhood, from
which it follows that $B_\m(\mathcal{S},\eta)$ is Stein. After
embedding it in some $\cx^N$ and using a retraction onto the
embedded submanifold, this reduces to \cite{wermer} Theorem 17.1.
(It is known that the smoothness assumed in this result is not the
best possible.) In our application, we will only be require the case
in which  $\mathcal{S}$ is diffeomorphic to the circle.

Now we turn to the proof of Proposition~\ref{smootharc}. We will in
fact prove the following proposition:

\begin{prop}\label{smoothnbds}
Let $k\geq 2$ and let $\alpha:S^1\rightarrow\m$ be a $\mathcal{C}^k$
embedding of the circle . Then the image $\alpha(S^1)$ has a
coordinate neighborhood  $W$ in $\m$ such that there is a coordinate
map $(\phi_1,\ldots,\phi_n):W\rightarrow\cx^n$ with each
$\phi_j\circ\alpha$  a $\mathcal{C}^{k}$ embedding of $S^1$ into
$\cx$.
\end{prop}

Indeed, any embedding of the interval can be extended to an
embedding of the circle, so Proposition~\ref{smoothnbds} immediately
implies Proposition~\ref{smootharc}.

\begin{proof} It is sufficient to consider
the case of $k=2$. Denote by $A_\delta$ the $\delta$ neighborhood
$B_{\cx^n}(S^1\times 0_{\cx^{n-1}}) $ of the circle
$S^1\times0_{\cx^{n-1}} $ in $\cx^n$. For small $\delta>0$ we will
construct a biholomorphic map $\Phi_\delta$ from $A_\delta$ onto an
open subset of $\m$ such that the image of $\Phi_\delta$ will
contain the embedded circle $\alpha(S^1)$. Consequently
$\Phi_\delta^{-1}$ is a coordinate map in a neighborhood of
$\alpha(S^1)$.

For $\eta>0$, let $\mathcal{X}_\eta= B_\m(\alpha(S^1),\eta)$. Also,
for a vector field $V$ on a manifold and a point $p$ on the manifold
let $\exp_Vp$ be the point to which $p$ flows in unit time along the
field $V$, that is $X(1)$, where $X(0)=p$ and $X'(t)= V(X(t))$. The
map $\exp_V p$ depends holomorphically on the vector field $V$ and
the point $p$.

We define a map from $A_\delta\subset\cx^n$ to
$\mathcal{X}_\eta\subset\m$ by setting:
\[ \Phi_\delta(z_1,\ldots,z_n)=\exp_{\sum_{j=2}^n
z_jf_j}\alpha_\delta(z_1), \]

where the number $\eta>0$, the vector fields $\{f_j\}_{j=2}^n$ on
the open submanifold $\mathcal{X}_\eta$,  and the map
$\alpha_\delta:B_\cx(S^1,\delta)\rightarrow\mathcal{X}_\eta$ are as
follows:

 (1) The holomorphic vector fields $\{f_j\}_{j=2}^n$, are
such that for each $z\in S^1\subset \cx$, the set of vectors
$\{\alpha'(z),f_2(\alpha(z)),\ldots,f_n(\alpha(z))\}$ spans the
tangent space $T_\alpha(z)\m$ over $\cx$.

To see that such $f_j$ exist, we note that $\mathcal{X}_\eta$ is
diffeomorphic to an open solid torus in $\cx^n$, and is Stein for
small $\eta$. Therefore by an application of the Oka principle, the
tangent bundle $T\mathcal{X}_\eta$ is trivial. Also, thanks to
Lemma~\ref{totallyreal}, any continuous function on the one
dimensional totally real submanifold $\alpha(S^1)$ may be
approximated by holomorphic functions in some neighborhood.
Therefore, the existence of the $f_j$ follows on approximating
smooth vector fields $\{g_j\}_{j=2}^n$ on $\alpha(S^1)$ such that
the set $\{\alpha'(z),g_2(\alpha(z)),\ldots,g_n(\alpha(z))\}$ spans
$T_{\alpha(z)}\m$, and shrinking $\eta$ to ensure the holomorphic
approximants $f_j$ are defined on $\mathcal{X}_\eta$.

    (2) Now we specify the map $\alpha_\delta$. This will be a
holomorphic map defined on $B_\cx(S^1,\delta)$ and taking values in
$\mathcal{X}_\eta\subset\m$ such that for some $0<\theta<1$:

\begin{itemize}
    \item on $S^1$ we have $\tmop{dist}(\alpha_\delta,\alpha)<C\delta^\frac{3}{2}$,
    as well as$\tmop{dist}(\nabla\alpha_\delta,\nabla\alpha) <C\delta^\frac{1}{2}$.
    \item there is a  constant $C$ (independent of $\delta$ ) such that and on
    $B_\cx(S^1,\delta)$ we have $\norm{\alpha_\delta}_{\mathcal{C}^{1,\theta}}<C$.
\end{itemize}

To construct $\alpha_\delta$ we note that since $\mathcal{X}_\eta$
is Stein, there is an embedding $j:\mathcal{X}_\eta\rightarrow\cx^M$
for large $M$, and there is a holomorphic retraction of a
neighborhood of $j(\mathcal{X}_\eta)$ onto $\mathcal{X}_\eta$.Fix
$\theta$, where $0<\theta<1$. Since $\alpha$ is of class
$\mathcal{C}^2$, we can use Lemma~\ref{analyticapprox} above to find
a holomorphic approximation $\alpha_\delta$ defined in a $\delta$
neighborhood of the circle $B_\cx(S^1,\delta)$ of $S^1$ in $\cx$,
and taking values in $\mathcal{X}_\eta$ such that the two conditions
above are satisfied.

For small $\delta$, the set
$\{\alpha_\delta'(z),f_2(\alpha_\delta(z)),\ldots,f_n(\alpha_\delta(z))\}$
spans $T_{\alpha_\delta(z)}\m$ for $z\in B_\cx(S^1,\delta)$.
Moreover, for small $\delta$, the map $\alpha_\delta$ is an
embedding. Therefore, for small enough $\delta$ the map
$\Phi_\delta$ is well defined and is a biholomorphism from
$A_\delta$ into $\mathcal{X}_\eta$.  Since the
$\mathcal{C}^{1,\theta}$ norm of $\alpha_\delta$ on
$B_\cx(S^1,\delta)$ is bounded independently of $\delta$, we
conclude that  $\Phi_\delta$ must be bounded in the
$\mathcal{C}^{1,\theta}$ norm on $A_\delta$. Recall that the tangent
bundle of $T\mathcal{X}_\eta$ is holomorphically trivial, and fix a
trivialization. Then $\Phi_\delta':A_\delta\rightarrow
\tmop{Mat}_{n\times n}(\cx)$ is a $\mathcal{C}^\theta$ map.
 Therefore for a constant $C_1$  independent of $\delta$
 and any $Z$ and $W$ in $A_\delta$  we have $\norm{\Phi_\delta'(W)-\Phi_\delta'(Z)}_{\mbox{op}} \leq
C_1\norm{W-Z}^\theta.$

In particular, if $Z$ lies on the circle $S^1\times 0_{\cx^{n-1}}$,
and $W$ is in the ball $B_{\cx^n}(Z,\delta)\subset A_\delta$, then
we will have
\begin{equation}\label{phideltaholder}
\norm{\Phi_\delta'(W)-\Phi_\delta'(Z)}_{\mbox{op}} \leq
C_1\delta^\theta.
\end{equation}

We claim that there is a constant $C_2$ independent of $\delta$ such
that if  $\delta>0$ is small, for every $Z\in A_\delta$ and every
tangent vector $v$ we have
\begin{equation}\label{phiprimeboundedbelow}
\norm{\Phi_\delta'(Z)v}\geq C_2\norm{v}.
\end{equation}

To see this, let $\tilde{\alpha}$ be an extension of $\alpha$ to a
neighborhood of $S^1$ in $\cx$ such that we have $\norm{\nabla
\tilde{\alpha}-\nabla\alpha_\delta}=O(\delta^\frac{1}{2})$ (see the
first conclusion of Lemma~\ref{analyticapprox}. The map
$\tilde{\alpha}$ is $\mathcal{C}^2$ and
$\overline{\partial}\tilde{\alpha}$ vanishes along $\Gamma$). We
define a map $\tilde{\Phi}$ by setting:
\[
\tilde{\Phi}(z_1,\ldots,z_n)=\exp_{\sum_{j=2}^n z_j
f_j}\tilde{\alpha}(z_1).
\]
Then $\tilde{\Phi}$ is a diffeomorphism from a neighborhood $A$ of
$S^1\times 0_{\cx^{n-1}}$ in $\cx^n$ into $\m$, and satisfies
$\norm{\Phi_\delta'-\tilde{\Phi}'}=O(\delta^\frac{1}{2})$. Clearly,
there is a constant $C>0$ such that for any $Z\in A$ and tangent
vector $v$ we have $\norm{\tilde{\Phi}'(Z)(v)}\geq C\norm{v}$. The
existence of the constant $C_2$ of
estimate~\ref{phiprimeboundedbelow} now follows immediately.

Ay an application of Lemma~\ref{babysurjective} above to the
inequalities it follows that  there is a $\delta_0>0$ , and a
constant $K$ independent of $\delta$ such that for
$\delta<\delta_0$, and $Z\in S^1\times 0_{\cx^{n-1}}$ we have
$\Phi_\delta(B_{\cx^n}(Z,\delta))\supset B_\m(\Phi_\delta(Z),
K\delta)$. Since for a point of the form $Z=(z,0,\ldots)\in
S^1\times 0_{\cx^{n-1}}$, we have $\Phi_\delta(Z)=\alpha_\delta(z)$,
we see that $\Phi_\delta(A_\delta)\supset
B_{\m}(\alpha_\delta(S^1),K \delta)$. On the other hand, for $z\in
S^1$,
\[
\mbox{dist}_{\m}\left(\alpha(z),\Phi_\delta(z,0,\ldots,0)\right)
\leq \mbox{dist}_{\m}\left(\alpha(z),\alpha_\delta(z)\right)
=O(\delta^{\frac{3}{2}}).
\]

Therefore for small $\delta$ we have, $ \alpha(S^1)\subset
\Phi_\delta(A_\delta)$, which shows that
$\Phi_\delta^{-1}:\Phi_\delta(A_\delta)\rightarrow\cx^n$ is a
coordinate map (biholomorphism onto an open subset of $\cx^n$)
defined in the neighborhood $\Phi_\delta(A_\delta)$ of $\alpha(S^1)$
in $\m$.

Note that as $\delta\rightarrow0$, the maps
$\Phi_\delta^{-1}\circ\alpha_\delta\rightarrow j$ on $S^1$ in the
$\mathcal{C}^1$ sense, where $j$ denotes the embedding of $S^1$ in
$\cx^n$ as $S^1\times 0^{n-1}$. Since
$\alpha_\delta\rightarrow\alpha$ in $\mathcal{C}^1$, it follows that
for small $\delta$, the first coordinate of
$\Phi_\delta^{-1}\circ\alpha$ is an embedding of $[0,1]$ into $\cx$.
Writing $\Phi_\delta^{-1}$ in coordinates as
$(\phi_1,\ldots,\phi_n)$, therefore, $\phi_1\circ\alpha$ is an
embedded arc in $\cx$ (which is obviously of class $\mathcal{C}^2$.
We now consider the coordinate system $(\phi_1, \phi_2 +
K\phi_1,\ldots, \phi_n + K\phi_1)$ in which for $j>1$, $\phi_j$ is
replaced by $\phi_j + K \phi_1$. For large $K$, every coordinate of
this map is a $\mathcal{C}^1$ embedding when restricted on $\alpha$.
\end{proof}

\subsection{Coordinate neighborhoods of $\mathcal{C}^1$ arcs}

We now prove Proposition~\ref{c1}. The first step is to show that
the image $\alpha([0,1])$ has a {\em Stein} neighborhood. We begin
with the following elementary observations : \begin{obs}
\label{spsh}
    Let $\gamma:[0,1]\rightarrow\m$ be an arc. Then $\gamma([0,1])$
    has a neighborhood $W$ such that there is a strictly plurisubharmonic
    function $\rho$ defined on $W$.
\end{obs}
Note that no regularity assumption apart from injectivity has been
made on $\gamma$.

\begin{proof}
There is of course a strictly \psh function in a neighborhood of
$\gamma(0)$. Suppose that for some $0<p<q<1$ the segment
$\gamma([p,q])$ is in a coordinate chart of $\m$, and $\rho$ is a
strictly \psh function in a neighborhood of $\gamma([0,p])$. By an
induction on a cover of $\gamma([0,1])$ by coordinate charts, it is
sufficient to construct a strictly \psh function in a neighborhood
of $\gamma([0,q])$. After subtracting a constant,
$\rho(\gamma(p))=0$, and there is a coordinate map $Z$ on a
neighborhood of $\gamma([p,q])$ so that $Z(p)=0\in\cx^n$. Fix an $r$
with $p<r<q$ so that $\rho$ is defined on $[p,r]$. Now we can define
the function $\tilde{\rho}$ as follows:

\[
\tilde{\rho} =\left\{ \begin{array}{cl}
\rho& \mbox{ near $\gamma([0,p])$}\\
\max\left(\rho,K\norm{Z}^2 -1\right) &\mbox{near $\gamma([p,r])$ with $K$ large (see below)}\\
K\norm{Z}^2 -1& \mbox{near $\gamma([r,q])$}
\end{array}\right.
\]
We  take $K$ in the above expression so large that
$K\norm{Z(\gamma(r))}^2 -1 > \rho(\gamma(r))$. Then $\tilde{\rho}$
is $\rho$ near $\gamma(p)$ and $ K\norm{Z}^2 -1$ near the other
endpoint $\gamma(r)$ of $\gamma([p,r])$ and continues to the next
chart.
\end{proof}
The following is a well known general fact regarding polynomially
convex sets:
\begin{obs}\label{polycvxpsh}
Let $X$ be a compact polynomially convex subset of $\cx^N$, an
$\lambda\geq 0$ be a continuous function on  $\cx^N$ such that
$\lambda=0$ exactly on $X$. Then, given any neighborhood $W$ of $X$,
there is a continuous plurisubharmonic function $\rho\geq 0$ such
that $\rho=0$ exactly on $X$, and on $W$, we have $\rho<\lambda$.
\end{obs}
\begin{proof}
For each $p\in \cx^n\setminus X$ there exists a continuous \psh
function $\rho_p\geq 0$ defined on $\cx^n$ such that $\rho_p(p)>0$,
and $\rho$ vanishes in a neighborhood of $X$. There exists a
sequence $x_j$ in $\cx^n\setminus X$ such that for all $z\in
\cx^n\setminus X$ there exists an integer $j$ such that
$\rho_{x_j}(z)>0$. e set $\rho = \sum \epsilon_j \rho_{x_j}$, with
$\epsilon_j>0$ small enough to ensure that $\epsilon_j \rho_{x_j}
\leq 2^{-j}\lambda$ on $W \cup B_{\cx^n}(0,j)$.
\end{proof}
We will now prove the following lemma:
\begin{lemma}\label{c1steinnbd}
Let $\gamma:[0,1]\rightarrow\m$ be a $\mathcal{C}^1$ embedded arc in
a complex manifold $\m$.  Then there is a $\delta>0$ with the
following properties:\\
(a) Let $0\leq t_0<t_1\leq 1$ be  such that $\abs{t_0-t_1}<\delta$,
and let $W_0$ and $W_1$ be given neighborhoods of $\gamma(t_0)$ and
$\gamma(t_1)$ in $\m$, Then there is a \psh function $\rho\geq 0$
defined in a neighborhood of $\gamma([t_0,t_1])$ in $\m$ such that
$\rho^{-1}(0)=V_0\cup \gamma([t_0,t_1])\cup V_1$, where $V_j\subset
W_j$ are compact neighborhoods of $\gamma(t_j)$ (for $j=0,1$).
\\
(b) Let $t_1$ be such that $0<t_1<\delta$ (resp. $0<1-t_1<\delta$),
and $W_1$ be a neighborhood of $\gamma(t_1)$ in $\m$. Then there is
a \psh function $\rho\geq 0$ defined in a neighborhood of
$\gamma([0,t_1])$ (resp. $\gamma([t_1,1])$) such that
$\rho^{-1}(0)=V_1\cup\gamma([0,t_1])$ (resp.
$\rho^{-1}(0)=\gamma([t_1,1])\cup V_1$. )
\end{lemma}
\begin{proof} Using compactness and local coordinates, it is clear
that is suffices to prove the result for $\m= \cx^n$. We only prove
statement (a) above, since the proof of (b) involves only minor
changes.

Let $t_0\in [0,1]$. After a linear change of coordinates, we can
assume that for each $j=1,\ldots,n$, the components of the tangent
vector $\gamma_j'(t_0)$ are non-zero. We let $\delta_{t_0}$ so small
that on $[t_0-\delta_{t_0},t_0+\delta_{t_0}]$ the component
functions $\gamma_j $ are $\mathcal{C}^1$ embeddings into $\cx$.

Now suppose we are given neighborhoods $W_0$ and $W_1$ of
$\gamma(t_0)$ and $\gamma(t_1)$ in $\m$. Choose $r>0$ so small that
for $j=1,\ldots,n$ the closed discs
$\overline{B_\cx(\gamma_j(t_0),r)}$ and
$\overline{B_\cx(\gamma_j(t_1),r)}$ are contained in the sets
$\pi_j(W_0)$ and $\pi_j(W_1)$ respectively, where
$\pi_j:\cx^n\rightarrow\cx$ is the $j$-th coordinate function. For
small $r$, the subset of $\cx$ 
\[ K_j= \overline{B_\cx(\gamma_j(t_0),r)}\cup\gamma_j([t_0,t_1])\cup
\overline{B_\cx(\gamma_j(t_1),r)}\] is polynomially convex, and
therefore there is a \psh $\rho_j\geq 0$  on $\cx$ which vanishes
exactly on $K_j$. Let $\rho \assign \sum_{j=1}^n \rho_j\circ\pi_j$.
Then clearly $\rho^{-1}(0)$ is the union of the subarc
$\gamma([t_0,t_1])$ with two closed polydiscs of polyradius $r$
centered at the endpoints $\gamma(t_0)$ and $\gamma(t_1)$, which are
contained in $W_0$ and $W_1$ respectively. Choosing $\delta$
uniformly for all $t_0$ by compactness, conclusion (a) follows.
\end{proof}
We can now prove Proposition~\ref{c1}.

It is clear that $\alpha$ is an embedded arc.  We consider two
partitions of the interval $[0,1]$
\[ 0<t_1<\cdots<t_{N-1}<1, \]
and
\[ 0<t_1'<\cdots<t_{N-1}'<1,\]
such that for each $j$ we have $t_j\not = t_j'$. We set $t_0=t_0'=0$
and $t_N=t_N'=1$. We will choose the partitions in such a way that
$\abs{t_{j}-t_{j+1}}<\delta$ and $\abs{t_j'-t_{j+1}'}<\delta$ for
$j=0,1,\ldots,N-1$, where $\delta$ is as in the conclusion of
Lemma~\ref{c1steinnbd} above. Suppose that for $j=1,\ldots,N-1$, the
open neighborhoods $W_j$ and $W_j'$ of $\gamma(t_j)$ and
$\gamma(t_{j+1})$ are such that $W_j\cap W_j'=\emptyset$.

We now apply lemma~\ref{c1steinnbd}. For $j=1,\ldots,N-1$, let
$\rho_j$ be a nonnegative \psh function in a neighborhood of
$\alpha([t_j,t_{j+1}])$  which vanishes exactly on
$\alpha([t_j,t_{j+1}])\cup V_j\cup V_{j+1}$, where $V_j \subset W_j$
and $V_{j+1}\subset W_{j+1}$ contain the points $\alpha(t_j)$  and
$\alpha(t_{j+1})$ respectively. Let $\rho_0$ and $\rho_N$ be
non-negative and \psh on neighborhoods of $\alpha([0,t_1])$ and
$\alpha([t_N,1])$ respectively such that
$\rho_0^{-1}(0)=\alpha([0,t_1])\cup V_1$ and
$\rho_N^{-1}(0)=\alpha([t_{N-1},1])\cup V_{N-1}$, where
$\alpha(t_1)\in V_1\subset W_1$ and $\alpha(t_{N-1})\in
V_{N-1}\subset W_{N-1}$. There is a function $\rho\geq 0$ in a
neighborhood of $\alpha([0,1])$ which  to be equal to $\rho_j$ in a
neighborhood of $\alpha([t_j,t_{j+1}])$. Since $\rho$ is locally the
max of \psh functions, it is itself \psh and vanishes exactly on the
arc $\alpha([0,1])$ and on small neighborhoods of $\alpha(t_j)$
 contained in $W_j$( where $j=1,\ldots, N-1$.)

 The same way we obtain a \psh $\rho'\geq 0$ which vanishes
exactly  on the arc $\alpha([0,1])$ and small neighborhoods of
$\alpha(t_j')$
 contained in $W_j'$( where $j=1,\ldots, N-1$.) Then $\tilde{\rho}=\rho+\rho'$ is
 a \psh function in a neighborhood of $\alpha([0,1])$ which vanishes
 exactly on $\alpha([0,1])$. Let $\epsilon>0$ be small, and $\psi$
 be a strictly \psh function in a neighborhood of $\alpha([0,1])$.
 Then the open set $\Omega=\{\tilde{\rho}<\epsilon\}$ supports the strictly
 \psh exhaustion function $(\epsilon-\tilde{\rho})^{-1}+\psi$, and
 is consequently Stein. If $\epsilon$ is small, the submersion
 $\phi$ is defined on $\Omega$.There is an embedding
$j:\Omega\hookrightarrow\cx^N$ for large enough $N$. Let
$\tilde{j}:\Omega\hookrightarrow\cx^{N+1}$ be the map
$\tilde{j}(z)\assign\left(j(z),\phi(z)\right)$, where $\phi$ is the
good submersion associated with the arc $\alpha$, whose existence is
assumed in the hypothesis. Then $\tilde{j}$ is again an embedding.
Let $\mathcal{X}\assign\tilde{j}(\Omega)$. Then,\\
 \noindent$\bullet$~ $\mathcal{X}$
is a complex submanifold of $\cx^{N+1}=\cx^N\times\cx$.\\
\noindent$\bullet$~$z_{N+1}:\mathcal{X}\rightarrow\cx$ is
submersion.\\
\noindent$\bullet$ Let $\tilde{\alpha}= \tilde{j}\circ\alpha$. Then
$\tilde{\alpha}$ is a $\mathcal{C}^1$ embedded arc  in
$\mathcal{X}\subset\cx^{N+1}$, such that the last coordinate
$\alpha_{N+1}:[0,1]\rightarrow\cx$ is a $\mathcal{C}^1$ embedding,
and\\
\noindent$\bullet$
 $\phi\circ\alpha:[0,1]\rightarrow\cx$ is a $\mathcal{C}^1$ embedding.
Set $\Gamma=\phi(\alpha([0,1]))$ and let $\psi:\Gamma\rightarrow
[0,1]$ be the inverse $\psi=(\phi\circ\alpha)^{-1}$. We let
\begin{eqnarray*}
\tilde{\beta}(z)&\assign& (\beta(z),z)\\
&\assign & (j\circ\alpha\circ\psi(z),z)
\end{eqnarray*}

To prove our result, it is sufficient to show that
$\tilde{\beta}(\Gamma)$ has a neighborhood $W$ in $\mathcal{X}$ such
that $W$ is biholomorphic to an open subset of $\cx^n$ and there is
a biholomorphism $w=(w_1,\ldots,w_n)$ from $W$ into $\cx^n$ such
that $w_n=z_{N+1}|_{\mathcal{X}}$ (where $(z_1,\ldots,z_N)$ are the
coordinates of $\cx^{N+1}$ in which $\mathcal{X}$ is embedded).

We will construct the map $w$ by first defining it on a neighborhood
of $\tilde{\beta}(\Gamma)$ in $\cx^{N+1}$ and its restriction to
$\tilde{X}$ will provide us with the required biholomorphic map.

To do this let $\{g_i\}_{i=1}^{n-1}$ be smooth maps from $\Gamma$
into $\cx^N$ such that for each $z\in\Gamma$ they, along with
$\beta'(z)$ span the tangent space
$T_{\beta(z)}(j(\mathcal{X}))\subset\cx^N$. If the $\cx^{N+1}$
valued maps $\{f_i\}_{i=1}^{n-1}$ are formed from $g_i$ by taking
the last coordinate to be $0$, then the $f_i(z)$ along with the
vector $\tilde{\beta}'(z) = (\beta'(z),1)$ span the tangent space
$T_{\tilde{\beta}(z)}(\mathcal{X})$.

 Let $A$ be a $\tmop{Mat}_{n\times N}(\cx)$ valued smooth map
 on $\Gamma$ such that $A(z)g_i(z)= e_i$
for each $i=1,\ldots,n-1$. We can approximate $A$ uniformly on
$\Gamma$ by a holomorphic matrix valued map $B$ defined in a
neighborhood of $\Gamma$ in $\cx$. Now we consider the map $\Lambda$
\[ \Lambda(z_1,\ldots,z_{N+1})\assign\left(B(z_{N+1})
\left(\begin{array}{c}z_1\\\vdots\\z_N\end{array} \right),z_{N+1}
\right) \] which is defined in a neighborhood of the arc
 $\tilde{\beta}(\Gamma)$ in $\cx^{N+1}$. Its derivative
is given by the matrix
\[ \Lambda'(z_1,\ldots,z_{N+1})=\left( \begin{array}{cc} B(z_{N+1}),
& B'(z_{N+1})\left(\begin{array}{c}z_1\\\vdots\\z_N\end{array}\right)\\
    0,& 1 \end{array} \right) \]

By construction, this is surjective from
$T_{\tilde{\beta}(z)}\mathcal{X}\subset\cx^{N+1}$ to $\cx^n$ at each
point of $\tilde{\beta}(\Gamma)$, if the approximation $B$
 is close enough. Moreover, it is clearly continuous. Therefore,
 $\Lambda$  maps a neighborhood  of the arc $\tilde{\beta}$
 in $\mathcal{X}$ to $\cx^n$
biholomorphically, and its last coordinate is $z_{N+1}$. This
completes the proof. \hfill$\square$

\subsection{  Mildly Singular Arcs, Step I: Stein neighborhoods}
\label{mildsingstein} The remaining part of this section is devoted
to a proof of Theorem~\ref{singularities}. The proof is in several
steps. In this step  we establish the existence of certain Stein
neighborhoods $\Omega_\delta$ of the arc $\alpha([0,1])$. This
allows us to solve $\dbar$ equations in these neighborhoods. In the
next step ( subsection \ref{immersiongluing})  we establish a result
regarding the gluing together of immersions defined in neighborhoods
of compact sets $K_1$ and $K_2$ in a manifold to a single immersion
defined in a neighborhood of their union.  We use these two results
in subsection~\ref{endproof}  to obtain a proof of
Theorem~\ref{singularities}.

 The main result of this section is the
following:
\begin{lemma}
\label{newstein} Let $\alpha:[0,1]\rightarrow\m$ be a
$\mathcal{C}^2$ arc with mild singularities. Let $P\subset[0,1]$ be
the set of points where $\alpha$ is not smooth  and let $\phi$ be
the good submersion associated with $\alpha$. Let $\mathcal{U}$ be a
fixed neighborhood of $\alpha(P)$ in $\m$. For $\delta>0$
sufficiently small
 there is a neighborhood $\Omega_\delta$ of $\alpha([0,1])$ in $\m$ such that
\begin{itemize}
\item $\Omega_\delta$ is Stein,
\item  $\Omega_\delta$ contains the $\delta$-neighborhood of the arc
$\alpha([0,1])$, and
\item away from the nonsmooth points $\alpha(P)$ of $\alpha$, the set
$\Omega_\delta$ coincides with the $\delta$
neighborhood of the arc. More precisely,
$\Omega_\delta \subset \mathcal{U} \cup B_\m\left(\alpha([0,1]),\delta \right).$
\end{itemize}
\end{lemma}

We will assume (without any loss of generality) that the points $0$
and $1$ are in $P$. We will need to use the following lemma which
gives a simple condition for the union of two polynomially convex
sets to be polynomially convex. A proof may be found in
(\cite{stout}, p. 386, Lemma 29.21(a)). $\widehat{X}$ denotes the
polynomial hull of a compact set $X\subset \cx^N$.
\begin{lemma}\label{stoutlemma}
Let $X_1$ and $X_2$ be compact polynomially convex sets in $\cx^n$
and let $p$ be a polynomial such that $\widehat{p(X_1)}\cap
\widehat{p(X_2)} \subset\{0\}$. If $p^{-1}(0)\cap \left(X_1\cup
X_2\right)$ is polynomially convex, then $X_1\cup X_2$ is again
polynomially convex.
\end{lemma}
We will also require the following
\begin{obs}\label{patchingobs}{\rm
Let $\gamma:[0,1]\rightarrow \m$ be an arc, and let $0<s<t<1$.
Suppose that in a neighborhood of $\gamma([0,t])$ is defined a
strictly plurisubharmonic function $\nu\geq 0$ which vanishes
precisely on the arc $\gamma$, and on a neighborhood of
$\gamma([s,1])$ is defined a continuous plurisubharmonic $\mu\geq0$
which also vanishes precisely on $\gamma$. Further suppose that
where both $\mu$ and $\nu$ are defined, $\mu<\nu$. Then there is a
continuous plurisubharmonic $\lambda$ in a neighborhood of
$\gamma([0,1])$ such that $\lambda$ coincides with $\mu$ near
$\gamma([0,s])$, with $\nu$ near $\gamma([t,1])$, and is bounded
above by $\nu$ near $\gamma([s,t])$. }\end{obs}
\begin{proof}
Let $\psi\leq0$ be a function of small $\mathcal{C}^2$ norm such
that \begin{itemize}
\item $\nu+\psi$ is still plurisubharmonic.
\item $\psi$ is 0 except in a small neighborhood of $\gamma(t)$,
where it is negative.
\end{itemize}
We set
\[
\lambda=\left\{
\begin{array}{ccc}
\nu & {\mbox{ near} } & \gamma([1,s])\\
\max(\mu, \nu +\psi) &{\mbox{ near }} & \gamma([s,t])\\
\mu &{\mbox{ near }} & \gamma([t,1])
\end{array}
 \right.\]
This will be plurisubharmonic provided the definition makes sense.
Now near $\gamma(s)$, we have $\lambda=\nu$ since $\mu<\nu$, so that
$\lambda$ is continuous in a neighborhood of $\gamma([1,t])$. Near
$\gamma(t)$, we have $\nu+\psi<0$, consequently $\lambda=\mu$ there,
and therefore $\lambda$ defines a continuous function in a
neighborhood of $\alpha([0,1])$.
\end{proof}

Now we prove Proposition~\ref{newstein}. Let $p\in P$. In a
neighborhood of $q=\alpha(p)$ we can find a system of coordinates
$(z_1,\cdots,z_n)$ such that the last coordinate $z_n$ is equal to
$\phi$. Now consider a polydisc $W$ of the type:
\[ W= \{(z_1,\ldots,z_n)\colon \abs{z_j}< R \mbox{  for } j=1,\ldots,n-1 ; \abs{z_n}<r\}\]

where $r$ and $R$ are so chosen that
 \begin{itemize}
\item $r$ is much smaller than $R$,
\item $W\Subset V_q$, where $V_q$ is a polydisc centered at $q$ such
that $V_q\Subset \mathcal{U}$.
\item the arc $\alpha$ enters and
exits $\overline{W}$ exactly once transversally through the part of
the boundary $\partial W$ given by
\[ \{(z_1,\ldots,z_n)\colon \abs{z_j}< R \mbox{  for } j=1,\ldots,n-1 ;\abs{z_n}=r\}.\]
That such $r$ and $R$ exist follows easily from  the fact that
$\phi\circ\alpha=z_n\circ\alpha$ is smooth.
\end{itemize}

Consider the compact set $ K\assign\overline{W}\cup \left(
\alpha([0,1])\cap \overline{V_q}\right)$. This is the union of the
polydisc $\overline{W}$ with  two ``whiskers" (the two components of
$\alpha([0,1])\cap(\overline{V_q}\setminus{W})$).The projection of
$K$ on the last coordinate is of the form $\overline{B(0,r)}\cup
\phi(\alpha(I))$ for a subinterval of $I$ of $[0,1]$. By the choice
of $r$ and $R$ above, this set is the disc $\overline{B(0,r)}$
attached with two arcs, each at exactly one point. Two applications
of Lemma~\ref{stoutlemma} (with $p=z_n$) shows that $K$ is
polynomially convex. Applying Observation~\ref{polycvxpsh} to the
polynomially convex $K$ and the continuous function
$\tmop{dist}(.,K)^2$ we obtain a continuous \psh function $\mu\geq
0$ in a neighborhood of $\alpha(P)$  such that $\mu(z)\leq
\tmop{dist}_\m(z,\alpha([0,1])^2$. Moreover, $\mu=0$ exactly on a
disjoint union of ``Soup can with whiskers"-type of neighborhoods of
the points of $\alpha(P)$. We take $\nu$  to be the square of the
distance to $\alpha([0,1])$. Then, applying
Observation~\ref{patchingobs} to the \psh $\mu$ and strictly \psh
$\nu$ ( twice for each point of $\alpha(P)$)  we can obtain a
plurisubharmonic $\lambda$ such that $\lambda$ vanishes on the arc,
and away from $\alpha(P)$ $\lambda=\nu$ the square of  the distance.
We now define $\Omega_\delta=\{Q\in \m\colon \lambda(Q)<\delta^2\}$.
It is easily verified that $\Omega_\delta$ has the two geometric
properties required, i.e. it contains the $\delta$-neighborhood of
$\alpha([0,1])$ and is actually the $\delta$-neighborhood away from
$\alpha(P)$. To see that it is Stein, we note that for small
$\delta$, the set $\Omega_\delta$ has the strictly \psh exhaustion
$(\delta^2-\lambda)^{-1} +\rho$, where $\rho$ is a strictly \psh
function in a neighborhood of $\alpha([0,1])$. (See
Observation~\ref{spsh} above.)
\subsection{Mildly Singular Arcs, Step 2: Gluing of Immersions}
\label{immersiongluing}

We now prove the induction step which we will use to glue locally
defined coordinate maps to obtain a coordinate map in a neighborhood
of a mildly singular arc.

 Let $K_1$ and $K_2$ be given compact subsets
of $\m$. Suppose we are given immersions $\Phi$ and $\Psi$ from
neighborhoods of $K_1$ and $K_2$ respectively into $\cx^n$. The main
question considered in this step is whether there is an immersion
from a neighborhood of the union $K=K_1\cup K_2$ into $\cx^n$. In
the application, $K_1\cup K_2$ will be a mildly singular arc without
any singular points in $K_1\cap K_2$.
\subsubsection{Hypotheses on the sets $K_1$, $K_2$ and the maps
$\Phi$ and $\Psi$}\label{hypothesessection}
 Of course to conclude that the immersions can be glued we
need to add additional hypotheses. These hypotheses should
correspond to the intended application. The ones that we will use
are the following:\\

{\bf 1.~Intersection is a smooth arc}
The basic hypothesis is that the intersection $K_1\cap K_2$ should
be a smooth arc. More precisely,there is a $\mathcal{C}^3$ arc
$\alpha:[0,1]\rightarrow\m$ such that $K_1\cap K_2$ is its image.\\

{\bf 2.~Already glued in one coordinate (``Special") }  We will assume that the two
immersions have already been glued in one coordinate. More
precisely, suppose that $\Phi$ and $\Psi$ are the immersions from
neighborhoods of $K_1$ and $K_2$ into $\cx^n$, then we will assume
that the last coordinates $\Phi_n$ and $\Psi_n$ are equal in a
neighborhood of the arc $K_1\cap K_2$.\\

We will denote by $\phi$ the map from a neighborhood of $K$ into
$\cx$ which is equal to $\Phi_n$ near $K_1$ and $\Psi_n$ near $K_2$.
In order to simplify the writing, we will say that a map whose last
coordinate is $\phi$ is {\em special}. Therefore, the hypothesis is
that there $\Phi$ and $\Psi$ are special. We will insist that while
modifying $\Phi$ and $\Psi$ so that they become glued, the last
coordinate is always $\phi$, i.e. they are special.\\

{\bf 3.~Good submersion} We will assume that the map $\phi$ of the
last paragraph is a good submersion associated with the arc
$\alpha$. Further, the sets $\phi(K_2\setminus K_1)$ and
$\phi(K_1\setminus
K_2)$ are disjoint.\\

{\bf 4.~Ghost of (\ref{newstein})} This hypothesis will be
required in the last step of the proof. We assume that for small
$\delta>0$, and for a given relatively compact neighborhood
$\mathcal{U}$ of $(K_1\setminus K_2)\cup (K_2\setminus K_1)$, there
is a Stein open neighborhood $\Omega_\delta$ of $K$ which has the
following properties.
\begin{itemize}
\item $\Omega_\delta$ contains the $\delta$-neighborhood of $K$.
\item Near the arc $K_1\cap K_2$, the set $\Omega_\delta$ is in
fact the $\delta$ neighborhood. More precisely $\Omega_\delta\subset
\mathcal{U}\cup B_\m\left(K_1\cap K_2,\delta\right)$. Consequently,
there is a fixed compact $H$ independent of $\delta$ such that each
$\Omega_\delta\subset H$.
\end{itemize}
With these hypotheses we state the following proposition, whose
proof will be given in $\S$\ref{approxgluesection} and
$\S$\ref{cousinsection} below.

\begin{prop}\label{glueingtheorem}
There is a special immersion $\Xi$ from a neighborhood of $K$ into
$\cx^n$.
\end{prop}

\subsubsection{Step 1 of proof of Prop.~\ref{glueingtheorem} :
 Approximate gluing of special immersions}\label{approxgluesection}
 We fix the map $\Phi$ and for
each small $\delta>0$  modify the special immersion $\Psi$ to a new
special immersion $\Psi_\delta$ such that near the arc $K_1\cap K_2$
the difference $\Psi_\delta- \Phi$ is small.Once this ``approximate
solution" is obtained, in $\S$\ref{cousinsection} we solve a
standard Cousin problem to modify both $\Phi$ and $\Psi_\delta$ so
that they now match on the intersection, and the result is an
immersion.

 We state the goal of this section
as a proposition.
\begin{prop}\label{approximateglueing}
After possibly shrinking the sets $K_1$ and $K_2$ (in a way so that
their union  $K_1\cup K_2$ is always $K$), we can find a constant
$C>0$  such that for each $\delta>0$ small a there is a  special
immersion $\Psi_\delta$ into $\cx^n$ defined in a neighborhood of
$K_2$ which contains the $\delta$ neighborhood of $K_1\cap K_2$,
 such that $\norm{(\Psi_\delta')^{-1}}_{\tmop{op}}<C$ and
on $B_\m(K_1\cap K_2,\delta)$ we have $\norm{\Phi-\Psi_\delta}=
O(\delta^2)$.
\end{prop}
For convenience, we divide the proof into a number of steps.

\paragraph{Step 1 }

{\em In this step we shrink $K_1$ and $K_2$ and modify $\Psi$ to
$\tilde{\Psi}$ in such a way that the derivatives of the immersions
$\Phi$ and $\tilde{\Psi}$ match at one point. This leads to a new
transition function $\tilde{\chi}$ with nicer properties.}\\
 By hypothesis (3) of $\S$\ref{hypothesessection} above,
 $\phi\circ\alpha$ is a smooth embedded  arc in $\cx$ (where $\phi$ is
  the common last coordinate of
$\Phi$ and $\Psi$). It follows that there is a neighborhood $W$ of
$K_1\cap K_2=\alpha([0,1])$ on which both $\Psi$ and $\Phi$ are
injective, and therefore biholomorphisms onto the image of $W$.

After translating in $\cx^n$, we can assume that
$\Phi(\alpha(\frac{1}{2}))=0$ and $\Psi(\alpha(\frac{1}{2}))=0$. The
 transition map $\chi\assign\Phi\circ\Psi^{-1}$ from
the open set $\Psi(W)\subset \cx^n$ onto $\Phi(W)\subset \cx^n$ is
biholomorphic. Since each of $\Phi$ and $\Psi$ is special (that is,
each has last coordinate equal to $\phi$)  it follows that $\chi$
has the form $ \chi(Z,w)=(\xi(Z,w),w)$ with $w\in\cx$, and $Z,
\xi(Z,w)\in\cx^{n-1}$. Moreover, $0\in W\cap \chi(W)$, in fact
$\chi(0)=0$.

Set $A=\chi'(0)\in \tmop{Mat}_{n\times n}(\cx)$. Define $
\tilde{\Psi} \assign A\circ\Psi$. Then $\tilde{\Psi}$ is again a
special immersion from a neighborhood of $K_2$ into $\cx^n$, and its
restriction to $W$ is a biholomorphic map onto the image. We can
also define a new transition function $ \tilde{\chi}
\assign\Phi\circ\tilde{\Psi}^{-1} =\chi\circ A^{-1}$, which is a
biholomorphism from $\tilde{\Psi}(W)\subset\cx^n$ onto
$\Phi(W)\subset\cx^n$. This new $\tilde{\chi}$ has the same form as
$\chi$, that is
\begin{equation}\label{tildechiform}
\tilde{\chi}(Z,w)=(\tilde{\xi}(Z,w),w)
\end{equation}
where $w\in\cx$, and both of $Z$ and $ \tilde{\xi}(Z,w)$ are in
$\cx^{n-1}$. The additional feature (not present before) is that
$\tilde{\chi}'(0) =\id$.

The derivative of $\tilde{\chi}$ is given by
\[ \tilde{\chi}' =
\left( \begin{array}{cc}
\tilde{\xi}_Z& \tilde{\xi}_w\\
0 & 1
\end{array}\right)\]
where subscripts denote differentiation, with $\tilde{\xi}_Z\in
GL_{n-1}(\cx)$ and $\tilde{\xi}_w$ a vector of $n-1$ components.

Since, $\tilde{\chi}'(0) =\id$, we can shrink the compact sets $K_1$
and $K_2$ (while not changing their union $K$), and the neighborhood
$W$ of $K_1\cap K_2$, so that $\tilde{\chi}'\approx \id$ on
$\tilde{\Psi}(W)$, in the sense that there is a holomorphic
$v:\tilde{\Psi}(W)\rightarrow \tmop{Mat}_{n-1\times n-1}(\cx)$ such
that on $\tilde{\Psi}(W)$, we have that
\begin{equation}
\label{vdef} \tilde{\xi}_Z = \exp\circ v.
\end{equation}

\paragraph{Step 2 }
{\em In this step we obtain, for $\delta>0$ small, an approximation
of the map $\tilde{\xi}$ of equation~(\ref{tildechiform}) by a  by a
map $\hat{\xi}_\delta$ affine in every coordinate except the last
one, such that $\norm{\hat{\xi}_\delta - \tilde{\xi}}=O(\delta^2)$ near
 the arc $\alpha$.}\\
Let $\lambda=\tilde{\Psi}\circ\alpha$. Then
$\lambda:[0,1]\rightarrow\tilde{\Psi}(W)\subset \cx^n$ is a
$\mathcal{C}^3$ arc. The last coordinate $\lambda_n$ is the embedded
$\mathcal{C}^3$ arc $\phi\circ\alpha$. Denote its image
$\lambda_n([0,1])$ by $\Gamma$. Since $\lambda_n$ is an embedding,
we can define a $\mathcal{C}^3$ map
$\gamma:\Gamma\rightarrow\cx^{n-1}$ by setting
\[ \gamma(\lambda_n(t))\assign \left(\lambda_1(t),\ldots,\lambda_{n-1}(t)\right)\]

 Apply Lemma~\ref{analyticapprox} it to $\gamma$. Hence, for
$\delta>0$ small, we can find a holomorphic
$\gamma_\delta:B_\cx(\Gamma,\delta)\rightarrow\cx^{n-1}$
 from the $\delta$ neighborhood $B_\cx(\Gamma,\delta)$,
of $\Gamma$ in $\cx$ into $\cx^{n-1}$ such that $\gamma_\delta$ is
bounded in the $\mathcal{C}^2$ norm on $B_\cx(\Gamma,\delta)$, and
on $\Gamma$, we have $\norm{\gamma_\delta-\gamma} =
O(\delta^\frac{5}{2})$. We now define for small $\delta>0$ a
$\cx^{n-1}$ valued holomorphic map $\hat{\xi}_\delta$ on the open
set $\cx^{n-1}\times B_\cx(\Gamma,\delta)\subset \cx^n$ by setting,
\[
\hat{\xi}_\delta(Z,w)\assign \tilde{\xi}(\gamma_\delta(w),w)
+\tilde{\xi}_Z(\gamma_\delta(w),w)
\left(Z-\gamma_\delta(w)\right).\]

$\hat{\xi}_\delta$ is the first order Taylor polynomial of the
$\cx^{n-1}$ valued map $\tilde{\xi}(\cdot,w)$ of $n-1$ variables
around the point $\gamma_\delta(w)\in\cx^{n-1}$.

We can rewrite $\hat{\xi}_\delta$ as
\begin{equation}\label{hatxideltadef}
\hat{\xi}_\delta(Z,w) = f_\delta(w) + \exp g_\delta(w) Z
\end{equation}
where $f_\delta$ , $g_\delta$ are holomorphic maps defined on
$B_\cx(\Gamma,\delta)$. The $\tmop{Mat}_{(n-1)\times (n-1)}(\cx)$
valued map $g_\delta$ is given by
\begin{equation}\label{gdeltadef}
 g_\delta(w)\assign v(\gamma_\delta(w), w)
\end{equation}
where $v$ is as in equation~(\ref{vdef}) above, that is, $\xi_Z
=\exp\circ v$. The $\cx^{n-1}$ valued map $f_\delta$ is defined by
\begin{equation}\label{fdeltadef}
f_\delta(w) \assign \xi(\gamma_\delta(w),w)- \left(\exp
g_\delta(w)\right) \gamma_\delta(w).
\end{equation}
We now observe the following two  facts which we will be use later.\\
{\bf 1.~$f_\delta$ and $g_\delta$ are bounded in $\mathcal{C}^2$}
 Since on $B_\cx(\Gamma,\delta)$, the map $\gamma_\delta$ is bounded in the
$\mathcal{C}^2$ norm independently of $\delta$, the same will be
true of the functions $f_\delta$ and $g_\delta$. That is, there is a
constant $C>0$  independent of $\delta$ such that for $j=0,1,2$ we
have $\norm{f_\delta^{(j)}}<C$ and $\norm{g_\delta^{(j)}}<C$.\\
{\bf 2.~$\hat{\xi}_\delta - \tilde{\xi}$ is small} More precisely,
 {\em suppose that the point
${\mathbf{Z}}=(Z,w)$ is in the $\delta$ neighborhood
$B_{\cx^n}(\lambda([0,1]),\delta)$ of the arc $\lambda([0,1])=
\Psi(K_1\cap K_2)$ in $\cx^n$. Then we have
\begin{equation}
\label{firstapprox}
 \norm{\hat{\xi}_\delta({\mathbf Z})-\tilde{\xi}({\mathbf Z}) } =
 O(\delta^2).
\end{equation}
}

To see this, note that $(Z,w)\in B_{\cx^n}(\lambda([0,1]),\delta)$
means that there is a $t\in \Gamma$ such that
$\norm{Z-\gamma(t)}<\delta$, and $\abs{w-t}<\delta$. Therefore,
using the properties of $\gamma_\delta$, we have
\begin{eqnarray*}
\norm{Z-\gamma_\delta(w)}&\leq & \norm{Z-\gamma(t)}
+\norm{\gamma(t)-\gamma_\delta(t)}
 + \norm{\gamma_\delta'}_{\sup}\abs{t-w}\\
&\leq& \delta + O(\delta^\frac{5}{2}) + C\delta = O(\delta).
\end{eqnarray*}

Now applying Taylor's theorem to the first order Taylor polynomial
$\hat{\xi}_\delta(\cdot,w)$ of $\tilde{\xi}(\cdot,w)$ around the
point $\gamma_\delta(w)$, we see that \[
\norm{\hat{\xi}_\delta({Z},w)-\tilde{\xi}({ Z},w) }=
O(\norm{Z-\gamma_\delta(w)}^2) = O(\delta^2).\]
\paragraph{Step 3} {\em  We now construct an approximation
$\chi_\delta$ of $\tilde{\chi}$.}\\
 At this point we require to use a lemma regarding the
approximation of functions of one variable. In order not to
interrupt the flow of the proof we state it here, but postpone its
proof to  $\S$~\ref{proofofapproxincx}.
\begin{lemma}
\label{approxincx} Let $B_1$, $B_2$ and $B_3$ be compact subsets of
$\cx$ such that $B_1\cap B_3=\emptyset$ and $B_1\cap B_2$ is a
single point, which we call $z_0$.
 Let $0<\theta<1$ and let $p$ be a positive integer. Then if $L$ is a closed subset of
$B_1$ such that $L\cap B_2=\emptyset$ (that is, $z_0\not\in L$),
then there is a constant $C$ with the following property. For
$\delta>0$ small, if $f$ is a holomorphic function in the closed
$\delta$-neighborhood $\overline{B_\cx(B_1\cup B_2,\delta)}$ of
$B_1\cup B_2$ such that
\[ \norm{f}_{\mathcal{C}^{1,\theta}(\overline{B_\cx(B_1\cup B_2,\delta)})}\leq 1 ,\]
then there is a holomorphic $f_\delta$ defined in the $\delta$
neighborhood of $B\assign B_1\cup B_2\cup B_3$ such that
\[ \norm{f_\delta}_{\mathcal{C}^1(B_\cx(B,\delta))}\leq C \]
and on $B_\cx(L,\delta)$ we have
\[ \abs{f-f_\delta}< C\delta^p. \]
\end{lemma}

To apply Lemma~\ref{approxincx} to our situation, we will let $p=2$
, $ B_1 =\lambda_n([0,\frac{3}{4}])$, $ L=
\lambda_n([0,\frac{1}{2}])\subset B_1 $ and $ B_2 =
\lambda_n([\frac{3}{4},1]).$ Then as required, we have that $B_1\cap
B_2$ a single point $z_0=\lambda_n(\frac{3}{4})$, and $z_0\not\in
L$. Also, we have $B_1\cup B_2= \lambda_n([0,1])= \Gamma$. For $B_3$
we take a relatively compact neighborhood of
 $\phi(K_2\setminus K_1)$ such that $B_1\cap B_3=\emptyset$. Then $B=
B_1\cup B_2\cup B_3 \supset \phi(K_2)$.

Now the holomorphic functions $f_\delta$ and $g_\delta$ defined in
equations (\ref{fdeltadef}) and (\ref{gdeltadef}) above are
holomorphic in the closed delta neighborhood of $\Gamma= B_1\cup
B_2$. Moreover, thanks to the $\mathcal{C}^2$ boundedness of the
maps, for any $\theta$ with $0<\theta<1$ we actually have $
\norm{f_\delta}_{\mathcal{C}^{1,\theta}(\overline{B_\cx(B_1\cup
B_2,\delta)})}$ and $
\norm{g_\delta}_{\mathcal{C}^{1,\theta}(\overline{B_\cx(B_1\cup
B_2,\delta)})}$ bounded independently of $\delta$.

Therefore, by an application of Lemma~\ref{approxincx}
 we get holomorphic maps $F_\delta$ and $G_\delta$, defined on the $\delta$ neighborhood $A_\delta$
of $B= \Gamma\cup B_3$, such that$\{ F_\delta \}$ and $\{ G_\delta
\}$ are uniformly bounded in the $\mathcal{C}^2$ norm independent of
$\delta$, and  on the set $B_\cx(\lambda_n([0,\frac{1}{2}]),\delta)$
we have $\norm{F_\delta- f_\delta} = O(\delta^2)$ and
$\norm{G_\delta- g_\delta} = O(\delta^2).$

We shrink the sets $K_1$ and $K_2$ again so that $K_1\cap K_2
=\alpha([0,\frac{1}{2}])$. We now define the map $\xi_\delta$ from
$\cx^{n-1}\times A_\delta$ into $\cx^{n-1}$ by
\begin{equation}\label{xideltadef}
\xi_\delta(Z,w) \assign F_\delta(w) +\left( \exp G_\delta(w)\right)
Z,
\end{equation}
and let
\begin{equation}
\label{chideltadef} \chi_\delta(Z,w) \assign
\left(\xi_\delta(Z,w),w\right)
\end{equation}

We now observe the following properties of the maps $\chi_\delta$
and $\xi_\delta$\\
\noindent$\bullet$ $\xi_\delta(Z,w)$ (and consequently
$\chi_\delta(Z,w)$) is defined on the set
 $ \cx^{n-1}\times A_\delta$ which contains a  neighborhood
of $\tilde{\Psi}(K_2)$ in $\cx^n$.\\
\noindent$\bullet${\em $\chi_\delta$ is a biholomorphic automorphism
of the set $\cx^{n-1}\times A_\delta$, and there is a constant $C>0$
independent of $\delta$ so that
\begin{equation}\label{chideltainverse}
\norm{\left(\chi_\delta^{-1}\right)'}_{\tmop{op}}<C.
\end{equation}
}\\
The derivative of $\chi_\delta$ is of the form
\[
\chi_\delta'(Z,w) = \left(\begin{array}{cc} \exp G_\delta(w) & *\\
                    0 & 1
                \end{array} \right)
\]
showing that $\chi_\delta'$ is a local biholomorphism. Now suppose
that $Z'\in\cx^{n-1}$ and $w'\in A_\delta$. Then we can solve the
equation $\chi_\delta(Z,w)= (Z',w')$ explicitly to obtain the
representations $ Z = \exp(-G_\delta(w')) \left(Z'-
F_\delta(w')\right)$ and $w=w'$ This shows that $\chi_\delta$ is a
biholomorphism, and its inverse is given by
\[ \chi_\delta^{-1} (Z,w) = \left( \exp(-G_\delta(w)) \left(Z- F_\delta(w)\right), w\right).\]
By construction each of $F_\delta$ and $G_\delta$ is bounded in the
$\mathcal{C}^1$ norm, independently of $\delta$. The bound
(\ref{chideltainverse}) follows immediately.\\
\noindent$\bullet${\em In the $\delta$ neighborhood
$B_\cx(\tilde{\Psi}(K_1\cap K_2),\delta)$ of $\tilde{\Psi}(K_1\cap
K_2)$ we have that
\[
\norm{\tilde{\chi} -\chi_\delta} = O(\delta^2)
\]}
This follows immediately from the fact that the coefficients
$F_\delta$ and $\exp G_\delta$ of $\xi_\delta$ are $O(\delta^2)$
perturbations of the coefficients  $f_\delta$ and $\exp g_\delta$ of
$\hat{\xi}_\delta$.

\paragraph{Step 4}{\em End of the proof of Proposition~\ref{approximateglueing}}. We now define the promised map
\[ \Psi_\delta\assign \chi_\delta\circ \tilde{\Psi}. \]
We observe that
\begin{itemize}
\item
By construction $\chi_\delta$ is a biholomorphic map from the set
$\cx^{n-1}\times A_\delta$ into itself. Recall that $A_\delta=
B_\cx(\phi(K_1\cap K_2),\delta)\cup L$, where $L$ is a fixed compact
neighborhood of $\phi(K_2\setminus K_1)\subset \cx$, such that
$L\cap\phi(K_1\cap K_2)$ is a single point. It easily follows from
this that $\Psi_\delta$ is defined on a set of the form
$B_\m(K_1\cap K_2,\delta)\cup \mathcal{U}$, where $\mathcal{U}$ is a
fixed (independent of $\delta$ ) neighborhood of $K_2\setminus K_1$.
\item $\Psi_\delta^{-1}$ and $\tilde{\Psi}^{-1}$  exist locally, and satisfy the equation
$\Psi_\delta^{-1} = \tilde{\Psi}^{-1}\circ \chi_\delta^{-1}$. Using
the chain rule and using the fact that
$\norm{(\chi_\delta^{-1})'}_{\tmop{op}}\leq C$ for a $C$ independent
of $\delta$ proved above we see that there is a constant $C'$
independent of $\delta$ such that
$\norm{(\Psi_\delta^{-1})'}_{\tmop{op}}<C'$. Another application of
the chain rule gives us $\norm{(\Psi_\delta')^{-1}}_{\tmop{op}}<C'$.
\item On the $\delta$ neighborhood $B_\m(K_1\cap K_2,\delta)$ of $K_1\cap K_2$, we have
$ \norm{\Phi-\Psi_\delta}= O(\delta^2)$.
 \end{itemize}

This completes the proof of Proposition~\ref{approximateglueing},
except we still need to prove Lemma~\ref{approxincx}.
\hfill$\square$

\subsubsection{Proof of Lemma~\ref{approxincx}.}
\label{proofofapproxincx} We need now to prove
Lemma~\ref{approxincx}. We will require the following fact, a proof
of which can be found in  \cite{stein}( ) The notation
$\mathcal{C}^{k,\theta}(K)$ denotes the  space of functions on $K$
whose $k$-th order derivatives are H\"{o}lder continuous with
exponent $\theta$.

We let $C$ stand for any constant not depending on $\delta$. The
proof will be in two very similar steps, each of which will involve
the solution of a $\dbar$ problem in one variable.\\

\subparagraph{Step 1.}

The aim of the first step is to construct a holomorphic function
$g_\delta$ on the $\delta$-neighborhood of $B= B_1 \cup B_2\cup B_3$
with the following properties.\\
\noindent$\bullet$ The $\mathcal{C}^1$ norm of $g_\delta$ is bounded
independently of $\delta$, that is
$\norm{g_\delta}_{\mathcal{C}^1(B_\cx(B,\delta))}\leq C$, where $C$
is independent of $\delta$.\\
\noindent$\bullet$ By hypothesis that the intersection $B_1\cap B_2$
is a single point $z_0$. Denote henceforth the disc
$B_\cx(z_0,\delta)$ by $N_\delta$, and the holomorphic function
$f-g_\delta$ (defined on a neighborhood of $B_1\cup B_2$) by
$k_\delta$. Then
    $ \norm{k_\delta}_{\mathcal{C}^{1,\theta}(N_\delta)}
    =O(\delta^{p+2})$.\\

To construct $g_\delta$ let $\psi$ be a $\mathcal{C}^\infty_c$
cutoff on $\cx$ which is $1$ on a neighborhood of $B_1$ (therefore
in a neighborhood of the point $z_0$ ) and $0$ on a neighborhood of
$B_3$. We now set
\[ \lambda_\delta(z)\assign \frac{1}{\left(z-z_0\right)^{p+4}}\Db{\psi}f(z).\]
Observe that $\lambda_\delta$ is defined on  the closed $\delta$
neighborhood $\overline{B_\cx(B_1\cup B_2)}$ of $B_1\cup B_2$ (and
this justifies the subscript $\delta$ on $\lambda$.) After extending
by $0$ at points where $\psi=0$, we can assume that $\lambda_\delta$
is defined and smooth on $B_\cx(B,\delta)$. Moreover, since $f$ is
bounded by 1 in the $\mathcal{C}^{1,\theta}$ norm, it follows that
there is a constant $C$ (independent of $\delta$), so that
$\norm{\lambda_\delta}_{\mathcal{C}^{1,\theta}(\overline{B_\cx(B,\delta)})}<C$.

Thanks to Lemma~\ref{genunifextension}, there is a compactly
supported extension $\tilde{\lambda}_\delta$ of $\lambda_\delta$ to
$\cx$ such that
$\norm{\tilde{\lambda}_\delta}_{\mathcal{C}^{1,\theta}(\cx)}\leq C$
(with $C$  independent of $\delta$).

We now define $ g_\delta \assign  \psi f + (z-z_0)^{p+4} \left(
-\frac{1}{\pi z} * \tilde{\lambda}_\delta\right)$, where $\psi f$ is
assumed to be $0$ where $\psi=0$. Then\\
\noindent$\bullet$ $g_\delta$ is holomorphic in the
$\delta$-neighborhood of $B$,\\
\noindent$\bullet$ since
$\tilde{\lambda}_\delta$ is bounded in the
    $\mathcal{C}^{1,\theta}$ norm, it follows that
    $\norm{g_\delta}_{\mathcal{C}^{1,\theta}(B_\cx(B,\delta))}\leq
        C$, where $C$ is independent of $\delta$,\\
\noindent$\bullet$for small enough $t$ we have $k_\delta(z_0+t) =
t^{p+4}\left(-\frac{1}{\pi z} * \lambda \right).$
        The first factor $t^{p+4}$ is $O(\delta^{p+2})$ in the
        $\mathcal{C}^2$ norm on the disc $N_\delta$. The second term
        is bounded on this disc in the $\mathcal{C}^{1,\theta}$
        norm.
    Therefore, we easily have
    \begin{equation}\label{kdeltaprime}
    \norm{k_\delta}_{\mathcal{C}^{1,\theta}(N_\delta)}=
    O(\delta^{p+2}).
    \end{equation}

Therefore, the function $g_\delta$ and $k_\delta=f-g_\delta$ satisfy
the conditions stated in the beginning of this step.
 \subparagraph{Step 2.} Now we write $
 f= k_\delta + g_\delta$. Observe that $g_\delta$ is already defined in the
$\delta$-neighborhood of $B= B_1\cup B_2 \cup B_3$. To approximate
$f$ we will approximate $k_\delta=f-g_\delta$ (which is holomorphic
on the $\delta$-neighborhood of $B_1\cup B_2$) by a holomorphic
function $h_\delta$ defined in the $\delta$ neighborhood of $B $. We
will set $ f_\delta= h_\delta+ g_\delta.$ Then, for $f_\delta$ to
satisfy the conclusion of Lemma~\ref{approxincx} it is sufficient
that\\
\noindent$\bullet\norm{h_\delta}_{\mathcal{C}^{1}(B_\cx(B,\delta))}$
is bounded independently of $\delta$.\\
 \noindent$\bullet$ on the
set $B_\cx(L,\delta)$
 (where as in the hypothesis $L$ is a fixed subset of $B_1$ not containing the point $z_0$)
we have  $\abs{ h_\delta - k_\delta } = O(\delta^p).$

To construct $h_\delta$ we proceed as follows. For $\delta>0$ small,
there is a $\mathcal{C}^\infty$ function $\alpha_\delta$ defined on
the closed $\delta$ neighborhood of $B$ such that\\
\noindent$\bullet $~$0\leq \alpha_\delta \leq 1$, with
$\alpha_\delta\equiv 1$ on the
    neighborhood of the set $L$, and $\alpha_\delta\equiv 0$ on a fixed neighborhood  of
    $B_3$, and\\
\noindent$\bullet$ $\nabla \alpha_\delta$ is supported in $N_\delta=
B_\cx(z_0,\delta)$ and furthermore,
        $\alpha_\delta$ satisfies
    \begin{equation}\label{alphadeltabound}
        \norm{\alpha_\delta}_{\mathcal{C}^2(N_\delta)}=O\left(\frac{1}{\delta^2}\right).
    \end{equation}
Now define a smooth function $\mu_\delta$ on
$\overline{B_\cx(B,\delta)}$  by setting $\mu_\delta\assign
\Db{\alpha_\delta}\cdot k_\delta$, and extending by 0 outside
$N_\delta$. Thanks to equations (\ref{alphadeltabound}) and
(\ref{kdeltaprime}), it follows that $
\norm{\mu_\delta}_{\mathcal{C}^{1,\theta}(B_\cx(B,\delta))}
=O(\delta^p). $ An application of Lemma~\ref{genunifextension} leads
to the construction of a compactly supported extension
$\tilde{\mu}_\delta$ of $\mu_\delta$ to the whole of $\cx$ such that
\[ \norm{\tilde{\mu}_\delta}_{\mathcal{C}^{1,\theta}(\cx)}= O(\delta^p). \]
We can assume that the supports of the $\tilde{\mu}_\delta$'s lie in
a fixed compact of $\cx$ (independently of $\delta$).

Define a holomorphic $h_\delta$ on the  $\delta$ neighborhood
$B_\cx(B,\delta)$ by setting
\[ h_\delta \assign \alpha_\delta \cdot k_\delta + \left( -\frac{1}{\pi z}* \tilde{\mu}_\delta\right) \]
where $\alpha_\delta \cdot k_\delta$ is understood to be 0 if
$\alpha_\delta=0$ (even if $k_\delta$ is not defined). It is clear
that $\norm{h_\delta}_{\mathcal{C}^1(B_\cx(B,\delta))}\leq C$.

Now consider the $\delta$ neighborhood $B_\cx(L,\delta)$ of the set
$L$. For small $\delta$, we have $\alpha_\delta\equiv 1$ on this
set. Then on $B(L,\delta)$ we have
\[\abs{h_\delta-k_\delta}=\abs{\frac{1}{\pi z}* \tilde{\mu}_\delta}
=O(\delta^p).\]  This ends the proof of Lemma~\ref{approxincx}.
\hfill$\square$
 \subsubsection{End of the proof of
Prop.~\ref{glueingtheorem} :
 A Cousin problem on a neighborhood of $K_1\cup K_2$}\label{cousinsection}
\begin{lemma}\label{cousinproblem}
For $\delta>0$ small, there are holomorphic maps $H_1^\delta$ and
$H_2^\delta$ from neighborhoods of $K_1$ and $K_2$ respectively into
$\cx^n$, such that
\begin{itemize}
\item in a neighborhood of $K_1\cap K_2$, we have
$ \Phi+ H_1^\delta = \Psi_\delta + H_2^\delta$,
\item for $j=1, 2$ we have
$ \norm{\left(H_j^\delta\right)'}_{\tmop{op}} = O(\delta^\frac{1}{2})$, and
\item the last coordinates of $H_1^\delta$ and $H_2^\delta$ are  both $0$.
\end{itemize}
\end{lemma}
We note that this completes the proof of
Proposition~\ref{glueingtheorem}. Consider the map $\Xi$ defined in
a neighborhood of $K=K_1\cup K_2$ by
\[ \Xi =\left\{\begin{array}{cc}
\Phi + H_1^\delta & {\mbox{ near $K_1$}}\\
\Psi_\delta + H_2^\delta & {\mbox{ near $K_2$}}
\end{array}
\right.
\]
This is a well defined special holomorphic map. To show that it is
an immersion, it is sufficient to show that the derivative $\Xi'(Z)$
is an isomorphism of vector spaces for $Z$ near $K$. Near $K_2$ we
have $\Xi'= \Psi_\delta'\left(\id +
\left(\Psi_\delta'\right)^{-1}\circ\left(H_2^\delta\right)'
\right)$, but $\norm{(\Psi_\delta')^{-1}}_{\tmop{op}}\leq C $ (with
$C$ independent of $\delta$), and
$\norm{\left(H_2^\delta\right)'}_{\tmop{op}} =
O(\delta^\frac{1}{2})$, so for small $\delta$, the linear operator
$\Xi'$ is an isomorphism. The same conclusion holds in a
neighborhood of $K_1$ for $\Xi= \Phi+ H_1^\delta$.

The proof of \ref{cousinproblem}   will require the following two
lemmas:
\begin{lemma}
\label{l2}
    Let $\m$ be Stein and let $\Omega\Subset\m$ be open. There is
    a constant $C$ with the following property.
    Let $U\subset \Omega$, where $U$ is a Stein open subset
    of $\m$, and let $g$ be a smooth $\dbar$-closed $(0,1)$-form
    on $U$ which is in $L^2_{(0,1)}(U)$. Then there is a smooth
     $u:\omega\rightarrow\cx$ such that $\dbar u =g$ and $u$ satisfies the estimate
    \begin{equation}
        \norm{u}_{L^2(U)} \leq C \norm{g}_{L^2_{(0,1)}(U)}
    \end{equation}
\end{lemma}
The point of the lemma is that the constant $C$ does not depend on
the open set $U$, but only on the relatively compact $\Omega$. This
is well known in the case when $\m=\cx^n$ (see e.g.
\cite{hormander}.) The general case may be reduced to the Euclidean
case by embedding $\m$ in some $\cx^N$ (for details see
\cite{thesis}).
\begin{lemma}
\label{wermerineq}
    let $\Omega\Subset\m $ be a smoothly bounded domain. Then,
    there is a constant $C$ such that for
    any smooth function $w:\Omega\rightarrow\cx$, and any compact
    $K\subset\Omega$  we will have the inequality:
    \begin{equation}
    \label{wermereqn}
    \sup_K\abs{w}
     \leq {C}\left\{\frac{1}{{\tmop{dist}}(K,\Omega^c)^n} \norm{w}_{L^2(\Omega)}+
        {{\tmop{dist}}(K,\Omega^c)}\norm{\dbar w}_{L^\infty(\Omega)} \right\}.
    \end{equation}
\end{lemma}
\begin{proof}
After using local co-ordinates and scaling, we can see that it is
sufficient to establish the following inequality for smooth
functions $w$ defined on the closed unit ball $B$ in $\cx^n$:
\[
\abs{w(0)} \leq K\left\{ \norm{w}_{L^2(B)} + \sup_B
\left(\max_j\abs{\frac{\partial w}{\partial \overline{z}_j}}
\right)\right\}.
\]
This is proved in \cite{wermer}( p. 130, Lemma 16.7)
\end{proof}
Now we prove Proposition~\ref{cousinproblem}.

\begin{proof}
For convenience of notation, we will suppress $\delta$ whenever
possible, and write $\Phi =(\hat{\Phi},\phi)$, $\Psi_\delta=
(\hat{\Psi}_\delta,\phi)$, and $H_j^\delta= (h_j,0)$, where
$\hat{\Phi}, \hat{\Psi}_\delta$ and $h_j$ are all
$\cx^{n-1}$-valued. Then the result is equivalent to solving the
$\cx^{n-1}$ valued additive Cousin problem
\[ h_1 -h_2 = \hat{\Psi}_\delta - \hat{\Phi}\assign R_\delta \]
with  bounds on $h_1$ and $h_2$.  We do this using  the well known
standard method. It is clear that $\left(K_1\setminus K_2\right)\cap
\left(K_2\setminus K_1\right) =\emptyset$, and therefore, we can
find a smooth cutoff $\mu$ in a neighborhood of $K_1\cup K_2$ such
that $\mu$ is $1$ in a neighborhood of $K_1\setminus K_2$ and $0$ in
a neighborhood of $K_2\setminus K_1$. We get a smooth solution to
the Cousin problem given by
\[ \left\{ \begin{array}{cccc}
\tilde{h}_1&=& \mu R_\delta & \mbox{extended by 0 where $\mu=0$}\\
\tilde{h}_2 &=& (\mu -1) R_\delta &\mbox{extended by 0 where
$\mu=1$}
\end{array}
\right.
\]
Observe that $\norm{\tilde{h}_j} =O(\delta^2)$ for $j=1,2$.

We want a ``correction" $u$, defined in a neighborhood of $K_1\cup
K_2$, so that $h_j =\tilde{h}_j + u$ will be holomorphic, that is,
we want to solve the equations $\dbar(\tilde{h}_j +u ) =0$. Both of
these equations are equivalent to the $\dbar$ equation in a
neighborhood of $K_1\cup K_2$ given by
\begin{equation}\label{dbarcousin}
 \dbar u =g ,
\end{equation}
where $g$ is a $\cx^{n-1}$ valued $(0,1)$ form (defined below) on a
Stein neighborhood $\Omega_\delta$ of $K_1\cup K_2$ of the type
whose existence was assumed in the hypotheses, that is,
$\Omega_\delta$ contains a $\delta$-neighborhood of $K_1\cup K_2$,
and {\em is}  the $\delta$-neighborhood near $K_1\cap K_2$.  The
smooth form $g$ is defined by $
 g \assign -R_\delta \dbar \mu $
on the $\delta$ neighborhood of $K_1\cap K_2$, and is extended by
$0$ to $\Omega_\delta$. We therefore have
\begin{eqnarray*}
\norm{g}_{L^2_{(0,1)}(\Omega_\delta)}
&\leq& \norm{g}_{L^\infty} \left(\tmop{Vol}(\tmop{support}g)\right)^\frac{1}{2}\\
&=&O(\delta^2) (O(\delta^{2n-1}))^\frac{1}{2} = O(\delta^{n+
\frac{3}{2}}).
\end{eqnarray*}

Using Lemma~\ref{l2} we obtain a function $u$ on $\Omega_\delta$
such that for some $C$ independent of $\delta$,
\[
\norm{u}_{L^2(\Omega_\delta)}\leq  C
\norm{g}_{L^2_{(0,1)}(\Omega_\delta)} \leq C
\delta^{n+\frac{3}{2}}.\] For convenience, denote the $\delta$
neighborhood of $K=K_1\cup K_2$ by $K_\delta$, Then, by hypothesis,
$K_\delta\subset \Omega_\delta$, and it follows that
$\tmop{dist}\left(K_{\frac{\delta}{2}},\Omega_\delta^c\right)\geq\frac{\delta}{2}$.
We apply now inequality \ref{wermereqn} of Chapter~1 to $u$, to
conclude that
\begin{eqnarray*}
\norm{u}_{L^\infty(K_{\frac{\delta}{2}})} & \leq &
C\left\{\frac{1}{\tmop{dist}\left(K_{\frac{\delta}{2}},\Omega_\delta^c\right)^n}\norm{u}_{L^2(\Omega_\delta)}
+
\tmop{dist}\left(K_{\frac{\delta}{2}},\Omega_\delta^c\right)\norm{\dbar
u}_{L^\infty(\Omega_\delta)}
 \right\}\\
& \leq& C\left\{ \frac{1}{\delta^n}\cdot \delta^{n+\frac{3}{2}} +
\delta\cdot \delta^2\right\} \leq  C \delta^\frac{3}{2}.
\end{eqnarray*}
Now, define $h_j = \tilde{h}_j + u$. Then we have $\norm{h_j} \leq
C\delta^2+ C\delta^\frac{3}{2}$ Therefore, on $K_\frac{\delta}{2}$,
we have $\norm{h_j} = O(\delta^\frac{3}{2})$. Applying the Cauchy
estimate, the required bound on $h_j'$ (and therefore
$(H_j^\delta)'$) follows.
\end{proof}
\subsection{Mildly singular arcs, Step 3: End of proof of Theorem~\ref{singularities}}
\label{endproof}
 Let $\alpha:[0,1]\rightarrow\m$ be a mildly
singular arc, and let $\phi$ be the associated good submersion into
$\cx$. We begin by covering the compact set $\alpha([0,1])$ by a
finite cover of open sets $\{ U_i\}_{i=1}^N$ such that
\begin{enumerate}
\item on each $U_i$ is defined a coordinate map whose last coordinate is $\phi$, and
\item the parts of the arc $\alpha$ in the intersections $U_i\cap U_{i+1}$ are all smooth.
\end{enumerate}

A simple induction argument applied to this cover shows that it is
sufficient to consider the case when $\alpha([0,1])$ is covered by
charts $U$ and $V$, such that there are coordinates
$\Phi:U\rightarrow\cx^n$, $\Psi:V\rightarrow\cx^n$, each having last
coordinate $\phi$. On the intersection, $\alpha$ is $\mathcal{C}^3$,
and $\phi\circ\alpha:[0,1]\rightarrow\cx$ is a $\mathcal{C}^3$ arc.
Thanks to Theorem~\ref{newstein}, we have for small $\delta>0$,
Stein neighborhoods $\Omega_\delta$ exactly of the type required.
Therefore, we obtain an immersion $\Xi$ from a neighborhood of
$\alpha([0,1])$ into $\cx^n$, whose last coordinate is $\phi$. But
$\phi\circ\alpha$ is injective, so that $\Xi$ is also injective near
$\alpha([0,1])$, that is, there is a neighborhood of the arc on
which $\Xi$ is a  biholomorphism.\hfill$\square$
\section{Approximation of Maps into Complex Manifolds}
\label{approx}
 We will continue to denote by $\m$ a complex
manifold of complex dimension $n$ which has been endowed with a
Riemannian metric (as before, the actual choice of the metric will be
irrelevant.) In analogy with the notation of Section~\ref{ac} we
introduce the following conventions. For a compact $K$ in $\cx$, the
notation $\mathcal{H}(K,\m)$ denotes the space of holomorphic maps
from $K$ to $\m$. A map $f:K\rightarrow\m$ is in $\mathcal{H}(K,\m)$
iff there is an open set $U_f$ in $\cx$ with $K\subset U_f$ and a
holomorphic $F:U_f\rightarrow\m$, such that $F$ restricts to $f$ on
$K$. By $\mathcal{A}^k(K,\m)$ we denote the closed subspace of
$\mathcal{C}^k(K,\m)$ consisting of those maps which are holomorphic
in the topological interior $\tmop{int}K$ of the compact set $K$.
The space $\mathcal{A}^k(K,\m)$  will always be considered to have
the topology inherited from $\mathcal{C}^k(K,\m)$. When $\m$ is the
complex plane $\cx$ we will abbreviate $\mathcal{H}(K,\cx)$ and
$\mathcal{A}^k(K,\cx)$ by $\mathcal{H}(K)$ and $\mathcal{A}^k(K)$
respectively.

By a {\em  Jordan domain} $\Omega$ in the plane, we mean a domain
whose boundary $\partial\Omega$ consists of finitely many Jordan
curves (homeomorphic images of circles in the place).  A Jordan
domain is said to be {\em circular} if each component of
$\partial\Omega$ is a circle in the plane. A $\mathcal{C}^1$ domain
is a Jordan domain in which each component of $\partial\Omega$ is a
$\mathcal{C}^1$ embedded image of a circle.

We now state the approximation results in our new notation.
\begin{theorem}\label{ctsapprox}
Let $\Omega\Subset\cx$ be a Jordan domain. Then
$\hol{\overline{\Omega}}{\m}$ is dense in
$\mathcal{A}^0(\overline{\Omega},\m)$.
\end{theorem}

For an analogous result for $\mathcal{C}^k$ maps with $k\geq 1$,we
have to assume more regularity on the boundary:
\begin{theorem} \label{smoothapprox}
Let $\Omega\Subset\cx$ a $\mathcal{C}^1$ domain, i.e. it is bounded
by finitely many $\mathcal{C}^1$ Jordan curves. If $k\geq1$, the
space $\hol{\overline{\Omega}}{\m}$ is dense in
$\mathcal{A}^k(\overline{\Omega},\m)$.
\end{theorem}

Before we proceed to prove Theorems~\ref{ctsapprox}  and
\ref{smoothapprox}, we will show that the boundary regularity
required in ~\ref{ctsapprox} can be reduced. We will show that in
Theorem~\ref{ctsapprox} it is sufficient to consider the case when
$\Omega$ is a circular domain, i.e. it is
sufficient to prove the following:\\
{\bf Theorem~\ref{ctsapprox}$^\prime$}~{\em For a circular domain
$W$, the subspace $\hol{\overline{W}}{\m}$ is dense in
$\mathcal{A}^0(\overline{W},\m)$.}\\
We will require the following two facts from the theory of conformal
mapping : (a)~(K\"{o}be)~ Let $\Omega$ be a Jordan domain. Then
there is a circular domain which is conformally equivalent to
$\Omega$.(b)~(Carath\'{e}odory) Let $\Omega_1$ and $\Omega_2$ be
finitely connected Jordan domains, and
$f:\Omega_1\rightarrow\Omega_2$ a biholomorphism. Then $f$ extends
to a homeomorphism from $\Omega_1$ onto $\Omega_2$. (See, for
example, \cite{tsuji}, Theorems IX.35 and IX.2 respectively. In this
reference {\bf (b)}  is  stated for simply connected domains but the
proof readily extends to the multiply connected case)
\begin{lemma}
\label{circular} Theorem~\ref{ctsapprox}$^\prime$ and
Theorem~\ref{ctsapprox} are equivalent.
\end{lemma}
\begin{proof}
It is clear that Theorem~\ref{ctsapprox} implies
Theorem~\ref{ctsapprox}$^\prime$ , since  a circular domain is also
a Jordan domain. For the converse we proceed as follows. Let
$\Omega$ be a Jordan domain. Thanks to the two facts from the theory
of conformal mapping mentioned before this proposition, it follows
that there is a circular domain $W$ such that there is a
homeomorphism $\chi:\overline{\Omega}\rightarrow \overline{W}$ which
maps $\Omega$ conformally onto $W$. Let
$f\in\mathcal{A}^0(\overline{\Omega},\m)$. Then $f\circ\chi^{-1}$ is
in $\mathcal{A}^0(\overline{W},\m)$, so by hypothesis we can
approximate it uniformly by functions $g\in\hol{\overline{W}}{\m}$.
Since $\chi\in\cts{\overline{\Omega}}{\cx}$, thanks to a version of
Mergelyan's theorem (see \cite{greenekrantz}, Theorem 12.2.7.), it
can be approximated uniformly on $\overline{\Omega}$ by functions
$\widetilde{\chi}\in \hol{\overline{\Omega}}{\cx}$. Therefore,
$g\circ\widetilde{\chi}$ is a holomorphic map defined in a
neighborhood of $\overline{\Omega}$ which approximates $f$
uniformly. This establishes the lemma.
\end{proof}

\subsection{Approximation on Good Pairs}\label{goodpairsection}
This section is devoted to the development of some tools which will
be used in $\S$\ref{lastsection} (along with the results of
Section~\ref{coordinate}. We begin with a few definitions.

\begin{definition}\label{goodpairdefn}{\rm We say that a pair $(K_1,K_2)$ of
 compact subsets of $\cx$ is a {\em good pair} if the following hold:\\
{\bf(A)} $K_1$ and $K_2$ are ``well-glued" together in the sense
that
\[\overline{K_1\setminus K_2} \cap \overline{K_2\setminus
K_1}=\emptyset.\] {\bf(B)} $K_1\cap K_2$ has finitely many connected
components, each  of which  is star shaped.}
\end{definition}

We now state the basic approximation result which will be used in
the proof of Theorems~\ref{ctsapprox} and \ref{smoothapprox}.
\begin{theorem}\label{goodpairapprox}
Let $(K_1,K_2)$ be a good pair of compact sets, and let $ V$ be a
compact subset of $\cx$ disjoint from $K_1$ such that the following
holds. For a fixed  $k\geq 0$, given any
$g\in\mathcal{A}^k(K_2,\cx)$ and an $\eta>0$, there is a $g_\eta\in
\mathcal{A}^k({K_2\cup V},{\cx})$ such that
$\norm{g-g_\eta}_{\mathcal{C}^k(K_2)}<\eta$.

Let $f\in\mathcal{A}^k(K_1\cup K_2,\m)$ be such that each of the
sets $f(K_j)$, for $j=1,2$ is contained in a coordinate neighborhood
of $\m$. Then, given $\epsilon>0$ there is an $f_\epsilon\in
\mathcal{A}^k({K_1\cup K_2},{\m})$ such that
$\tmop{dist}_{\mathcal{C}^k({K_1\cup
K_2},{\m})}\left(f,f_\epsilon\right)<\epsilon$, and $f_\epsilon$
extends as a holomorphic map to a neighborhood of ${(K_2\cap V)}$.
\end{theorem}
Of course this  is of interest only in the case when $K_2\cap V\not
=\emptyset$. We split the proof into several steps.
\begin{obs}
\label{additive} (Additive Cousin problem $\mathcal{C}^k$ to the
boundary.) Let $(K_1,K_2)$ be a good pair. For each $k\geq 0$, there
exist bounded linear maps $T_j: \mathcal{A}^k(K_1\cap
K_2,\cx)\rightarrow \mathcal{A}^k(K_j,\cx)$ such that for any
function $f$ in $\mathcal{A}^k(K_1\cap K_2,\cx)$ we have on $K_1\cap
K_2$,
\begin{equation}
\label{additive_eq}
 T_1f +  T_2f = f,
\end{equation}
\end{obs}
\begin{proof} We reduce the problem to a $\dbar$
equation in the standard way. Let $\chi$ be a smooth cutoff which is
1 near $\overline{K_1\setminus K_2}$ and 0 near
$\overline{K_2\setminus K_1}$. Let $\lambda\assign f.\Db{\chi}$, so
that $\lambda\in \mathcal{A}^k(K_1\cap K_2,\cx)$. Let $0<\theta<1$.
Thanks to Lemma~\ref{genunifextension},
there is a bounded linear extension operator
$E:\mathcal{C}^k(K_1\cap K_2)\rightarrow\mathcal{C}^{k-1,\theta}(\rl^2)$, such
that each for $g$, the extension $Eg$ is supported in a fixed
compact set  of $\rl^2$.  We can now define  $ T_1 f = (1-\chi).f +
\frac{1}{\pi z}*(E\lambda)$, and  $T_2 f = \chi.f -\frac{1}{\pi
z}*(E\lambda)$, where $(1-\chi).f$ (resp. $\chi.f$) is assumed to be
$0$ at points where $\chi=0$ (resp. $(1-\chi)=0$) even if $f$ is not
defined. Since for $U\Subset\cx$, the map $v\mapsto\frac{1}{z}*v$ is bounded from 
$\mathcal{C}^{k-1,\theta}_c(U)$ to $\mathcal{C}^{k,\theta}(\cx)$,
 the result follows.
\end{proof}
 We use Observation~\ref{additive} to prove a  version of the Cartan
 Lemma on factoring matrices  similar to one  found in \cite{douady},
 pp. 47-48.
\begin{lemma}\label{cartan}
Let $(K_1, K_2)$ be a good pair, and let $g\in \mathcal{A}^k(K_1\cap
K_2,GL_n(\cx))$, where $k\geq 0$. Then, for $j=1,2$ there are
$g_j\in \mathcal{A}^k(K_j,GL_n(\cx))$ such that $g=g_2\cdot g_1$ on
$K_1\cap K_2$.
\end{lemma}
\begin{proof}
Denote by $\mathfrak{G}_j$ the  group $\mathcal{A}^k(K_j,GL_n(\cx))$
and by $\mathfrak{G}$  the group $\mathcal{A}^k(K_1\cap
K_2,GL_n(\cx))$. Let $\mu:\mathfrak{G}_1\times
\mathfrak{G}_2\rightarrow\mathfrak{G}$ be the map $\mu(g_1,g_2)=
g_2|_{K_1\cap K_2} \cdot g_1|_{K_1\cap K_2}$. Its derivative at the
point $(1_{\mathfrak{G}_1},1_{\mathfrak{G}_2})$ is given by the
linear map from the Banach space
$\mathcal{A}^k(K_1,\tmop{Mat}_{n\times
n}(\cx))\oplus\mathcal{A}^k(K_2,\tmop{Mat}_{n\times n}(\cx))$ into
the Banach space $\mathcal{A}^k(K_1\cap K_2,\tmop{Mat}_{n\times
n}(\cx))$ given by $(h_1,h_2)\mapsto h_1|_{K_1\cap K_2} +
h_2|_{K_1\cap K_2}$. Thanks to observation~\ref{additive} above,
this is surjective. Consequently, there is a neighborhood $U$ of the
identity in $\mathfrak{G}$ such that for any $g\in U$ there are
$g_j$ in $\mathfrak{G}_j$ such that $g=g_2 g_1$ on $K_1\cap K_2$.
This proves the assertion when $g$ is in the neighborhood $U$. We
may assume without any loss of generality that the exponential map
is surjective diffeomorphism onto $U$ from a neighborhood $V$ of $0$
in $\mathcal{A}^k(K_1\cap K_2,\tmop{Mat}_{n\times n}(\cx))$.

For the general case, observe that thanks to the fact that each
component of $K_1\cap K_2$ is contractible, the group $\mathfrak{G}$
is connected, and hence $\mathfrak{G}$ is generated by the
neighborhood $U$ of the identity. Therefore, we may write
$g=\prod_{i=1}^N\exp(h_i)$, where the $h_i\in V$. Now, the set
$\cx\setminus K_1\cap K_2$ is connected, and therefore it is
possible to approximate each $h_j$ by an entire matrix valued
$\tilde{h}_j$ such that on $K_1\cap K_2$. Let $\tilde{g} =\prod
\tilde{h}_j$, and $\hat{g}=\tilde{g}^{-1}\cdot g$. If the
approximation of $h_j$ by $\tilde{h}_j$ is close enough $\hat{g}\in
U$, and consequently it is possible to write $\hat{g}=a_2 a_1$,
where $a_j\in \mathfrak{G}_j$. We can take $g_1=a_1$ and
$g_2=\tilde{g}a_2$ to complete the proof.
\end{proof}
The following solution of a non-linear Cousin problem  is due to
Rosay (\cite{rosay1}, also see comments in \cite{rosay2}).
\begin{lemma}\label{rosay}
Let $\omega$ be an open subset of $\cx^n$ and let
$F:\omega\rightarrow\cx^n$ be a holomorphic immersion. Let
$(K_1,K_2)$ denote a good pair of compact subsets of $\cx$, and for
some $k\geq 0$, let $u_1\in\mathcal{A}^k(K_1,\cx^n)$ be such that
$u_1(K_1\cap K_2)\subset\omega$. Given any $\epsilon>0$, there
exists $\delta>0$ such that  if  $u_2\in\mathcal{A}^k(K_2,\cx^n)$ be
such that $\norm{u_2-F(u_1)}<\delta$, then for $j=1,2$  there exist
$v_j\in \mathcal{A}^k(K_j,\cx^n)$ such that $\norm{v_j}<\epsilon$,
and $u_2+v_2=F(u_1+v_1)$.
\end{lemma}
It is important to note that the map $u_1$ is {\em fixed}. In
\cite{rosay2}, a version is proved in which this restriction is
removed. This requires a version of Cartan's lemma for {\em bounded}
matrices (see (\cite{br})),a  result  valid  if $K_1$, $K_2$ and
$K_1\cup K_2$ are simply connected. Unfortunately, such a result
could not be proved for the
 more general  $K_1,K_2$ considered here.
 The proof will use the following well-known fact from the theory of
 Banach spaces, which can be proved using a standard iteration
 argument (see \cite{lang},  pp. 397-98)~:
\begin{lemma}
\label{surjective}
  Let $\mathcal{E}$ and $\mathcal{F}$ be Banach Spaces and let
$\Phi: B_{\mathcal{E}}(p,r)\rightarrow\mathcal{F}$ be a
$\mathcal{C}^1$ map. Suppose there is a constant $C>0$ such that:
\begin{itemize}
\item
for each $h\in B_{\mathcal{E}}(p,r)$, the linear operator
 $\Phi'(h):\mathcal{E}\rightarrow\mathcal{F}$ is surjective and the
equation $\Phi'(h)u=g$ can be solved for $u$ in $\mathcal{E}$ for
all $g$ in $\mathcal{F}$ with the estimate
$\norm{u}_{\mathcal{E}}\leq C \norm{g}_{\mathcal{F}}$.
\item
for any $h_1$ and $h_2$ in $B_{\mathcal{E}}(p,r)$ we have $\norm{
\Phi'(h_1)-\Phi'(h_2)}\leq \frac{1}{2C}.$
\end{itemize}
Then,
\[ \Phi(B_{\mathcal{E}}(p,r)) \supset B_\mathcal{F}\left( \Phi(p),\frac{r}{2C}\right).\]
\end{lemma}
We now give a proof of Lemma~\ref{rosay}.

 \begin{proof}
 Denote by $\mathcal{E}$ the Banach space
$\mathcal{A}^k(K_1,\cx^n)\oplus\mathcal{A}^k(K_2,\cx^n)$, which we
endow with the norm $\norm{\cdot}_{\mathcal{E}} \assign
\max\left(\norm{\cdot}_{\mathcal{A}^k(K_1,\cx^n)},\norm{\cdot}_{\mathcal{A}^k(K_2,\cx^n)}\right)$.
Also, let the open subset $\mathcal{U}$ of $\mathcal{E}$ be given by
$\{(w_1,w_2): w_1(K_1\cap K_2)\subset \omega\}$. Then
$(u_1,w_2)\in\mathcal{U}$, for any $w_2\in\mathcal{A}^k(K_2,\cx^n)$.
a Let $\mathcal{F}$ be the Banach space $\mathcal{A}^k(K_1\cap
K_2,\cx^n)$. Consider the map
$\Phi:\mathcal{U}\rightarrow\mathcal{F}$ given by
$\Phi(w_1,w_2)\assign w_2|_{K_1\cap K_2}- F\circ (w_1|_{K_1\cap
K_2})$.  A computation shows that $\Phi'(w_1,w_2)$ is the linear map
from $\mathcal{E}$ to $\mathcal{F}$ given by $(v_1,v_2)\mapsto
v_2|_{K_1\cap K_2}- F'(w_1|_{K_1\cap K_2})(v_1|_{K_1\cap K_2}).$
Observe that $w_2$ plays no role whatsoever in this expression, and
therefore  $\Phi'(w_1,w_2)\in BL(\mathcal{E},\mathcal{F})$ is in
fact a smooth function of $w_1$ alone, and we will henceforth denote
it by $\Phi'(w_1,*)$.

 We  construct a right inverse to $\Phi'(u_1,*)$ .
Let $\gamma=F'(u_1)|_{K_1\cap K_2}$. Then $\gamma\in
\mathcal{A}^k(K_1\cap K_2,GL_n(\cx))$, and thanks to
Lemma~\ref{cartan} above, we may write $\gamma= \gamma_1\cdot
\gamma_2$, where $\gamma_j\in \mathcal{A}^k(K_j,GL_n(\cx))$. (We
henceforth suppress the restriction signs.) For $g\in \mathcal{F}$,
let $S(g)=(-\gamma_1^{-1} T_1(\gamma_2^{-1}g),\gamma_2
T_2(\gamma_2^{-1}g))$, where the $T_j$ are as in
equation~\ref{additive_eq} above. Then $S$ is a bounded linear
operator from $\mathcal{F}$ to $\mathcal{E}$, and a computation
shows that $\Phi'(u_1,*)\circ S(g)$ is the identity map on
$\mathcal{F}$. Choose $\theta>0$ so small so that if
$w_1\in\mathcal{F}$ is such that $\norm{w_1 - u_1}<\theta$ then (a)
the equation $\Phi'(w_1,*)u=g$ can be solved with the estimate
$\norm{u}\leq 2 \norm{S} \norm{g}$, and (b) $\norm{\Phi'(w_1,*)-
\Phi'(u_1,*)}_{\tmop{op}}< \frac{1}{8\norm{S}}$. (These follow from
continuity and the fact that small perturbations of surjective
linear operator are still surjective) Consequently, if
$\epsilon<\theta$ and $u_2\in\mathcal{A}^k(K_2,\cx^n)$, for the ball
$B_{\mathcal{E}}((u_1,u_2),\epsilon)$ the hypothesis of
Lemma~\ref{surjective} are verified with $C=2\norm{S}$. We have
therefore,
\begin{eqnarray*}
\Phi\left(B_\mathcal{E}((u_1,u_2),\epsilon)\right)&\supset&
 B_\mathcal{F}\left(\Phi(u_1,u_2),\frac{\epsilon}{2C}\right)\\
& =&  B_\mathcal{F}\left(u_2- F(u_1),\frac{\epsilon}{2C}\right).
\end{eqnarray*}
So, if $\norm{u_2- F(u_1)}< \frac{\epsilon}{4C}$, we have
$0\in\Phi\left(B_\mathcal{E}((u_1,u_2),\epsilon)\right)$. This is
exactly the conclusion required.
 \end{proof}
 We will now prove Theorem~\ref{goodpairapprox}.
\begin{proof}
We omit the restriction signs on maps for notational clarity. For
$j=1,2$  let the coordinate neighborhoods $V_j$ of $\m$ be such that
$f(K_j)\subset V_j$. We begin by fixing biholomorphic maps
$F_j:V_j\rightarrow F_J(V_j)\subset\cx^n$, and set $F= F_2\circ
F_1^{-1}$. Then $F$ is a biholomorphism from the open set
$\omega=F_1(V_1\cap V_2)$ onto the open set $F_2(V_2\cap V_1)$.
Moreover, a pair of maps $w_1$ and $w_2$ from $K_1$ and $K_2$
respectively ``glue together" to form a map from $K_1\cup K_2$
(i.e., there is a map $h: K\rightarrow \m$ such that $w_j = F_j\circ
h$), only if $w_2 = F(w_1)$.

Let $u_1=F_1\circ f\in\mathcal{A}^k(K_1,\cx^n)$. Since $V\cap
K_1=\emptyset$ by hypothesis, the pair of compact sets $(K_1,
K_2\cup V)$ is good. Fix $\epsilon_0>0$ and  let $\delta_0>0$ be the
$\delta$ corresponding to $\epsilon=\epsilon_0$ in
Lemma~\ref{rosay}, for the good pair $(K_1,K_2\cup V)$ and
$F,\omega,u_1$ as defined above. Let $u_2\in \mathcal{A}^k(K_2\cup
V,{\cx^n})$ be a $\mathcal{C}^k$ approximation to
$F_2\circ\left(f|_{K_2}\right)$ such that $ u_2(K_2)\subset V_2$,
and  $\norm{u_2 - F(u_1)}<\delta_0$. (Such a $u_2$ exists by
hypothesis).

Then, by Proposition~\ref{rosay},  there is a
$v_1\in\mathcal{A}^k(K_1,\cx^n)$ and a $v_2\in\mathcal{A}^k({K_2\cup
V}{\cx^n})$, such that $\norm{v_j}<\epsilon_0$ and $u_2 + v_2 =
F(u_1 + v_1)$. Hence the maps $u_1+v_1$ and $u_2+v_2$ glue together
to form a map  $f_{\epsilon_0}$  given by
\[
f_{\epsilon_0}\assign \left\{ \begin{array}{ccc}
F_1^{-1}(u_1+v_1)& \tmop{ on } & K_1\\
F_2^{-1}(u_2+v_2) &  \tmop{ on } & K_2 \tmop{  and ~near } K_2\cap V
\end{array} \right.
\]
Clearly, $f_{\epsilon_0}$ is in $\mathcal{A}^k(K_1\cup K_2,\m)$, and
extends to a holomorphic map near $K_2\cap V$. Moreover,
$\tmop{dist}_{\mathcal{C}^k}(f_{\epsilon_0}, f)= O(\epsilon_0)$. The
result follows.
\end{proof}

\subsection{Proof of theorems~\ref{ctsapprox}$^\prime$ and
\ref{smoothapprox}}\label{lastsection}
 Let  $k\geq 0$, be an integer, and let the domain $\Omega$ be circular 
 if $k=0$, and $\mathcal{C}^1$ if $k\geq 1$. We fix 
$f\in\mathcal{A}^k(\overline{\Omega},\m)$. We wantto approximate
$f$ in the $\mathcal{C}^k$ sense on $\overline{\Omega}$.

The basic idea of this proof is to slice the  $\overline{\Omega}$ by
a system of parallel lines. If the slices are
narrow enough, we will show that thanks to the results of
Section~\ref{coordinate} the  graph of  $f$ over each slice 
is contained in a coordinate neighborhood of $\m$. We will further
show that the slicing can be done in a way
that the unions of alternate slices form a good pair. Consequently,
we can use the results of $\S$\ref{goodpairsection} to prove the
approximation results. We break the proof up into a sequence of
lemmas.

\begin{lemma}\label{slicing}Denote by $F$ the map in
$\mathcal{A}^k(\overline{\Omega},\cx\times\m)$ given by
$F(z)=(z,f(z))$. For real $\xi$, let  $L_\xi$ be the vertical
straight line $\{z\colon \Re z=\xi\}$. Then,
 \begin{itemize}\item $F(L_\xi\cap \overline{\Omega})$
has a coordinate neighborhood in $\cx\times\m$, and
\item There is a nowhere dense $E\subset\rl$  such that if
$\xi\not\in E$, the line $L_\xi$ meets $\partial\Omega$
transversely. If $k=0$, the set $E$ can even be taken to be finite.
\end{itemize}
\end{lemma}
\begin{proof} We first prove that $F(L_\xi\cap \overline{\Omega})$
has a coordinate neighborhood. Each connected component of
$L_\xi\cap \overline{\Omega}$ is a point or a compact interval, and
thanks to the injectivity of $F$, if we show that for each such
component $I$, the set $F(I)$ has a coordinate neighborhood in
$\cx\times\m$, it will follow that $F(L_\xi\cap \overline{\Omega})$
has a coordinate neighborhood in $\cx\times\m$. This is trivial if
the component $I$ is a point.
Therefore, let $I$ be a compact interval. We consider three cases:\\
{\em Case 1: $k=0$.}~ In this case $\Omega$ is a circular domain,
and hence for each $\xi$ the set $L_\xi\cap \partial{\Omega}$ and a
fortiori $I\cap \partial\Omega$ is finite. Consider the arc $F|_I$
in $\cx\times\m$. This arc is real analytic off the finite set of
points $I\cap \partial\Omega$, and the
projection~$\phi:\cx\times\m\rightarrow\cx$ has the property that
$\phi\circ\left(F|_I\right)$ is the inclusion map $I\hookrightarrow
\cx$. Consequently, $F|_I$ is a real analytic mildly singular arc in
the sense of Definition~\ref{mildsingularitiesdefn}, and hence
thanks to Theorem~\ref{singularities}, $F(I)$ has a coordinate
neighborhood.\\

{\em Case 2: $k=1$.}~ In this
case, the arc $F|_I$ is $\mathcal{C}^1$. As in the last case, let
$\phi$ be the projection $\phi:\cx\times\m\rightarrow\cx$, which has
the property that $\phi\circ\left(F|_I\right)$ is the inclusion map
$I\hookrightarrow \cx$. Therefore, Proposition~\ref{c1} applies,
and $F(I)$ has a coordinate neighborhood.\\

{\em Case 3: $k\geq 2 $.}~ In this case $F|_I$ is a $\mathcal{C}^k$
embedded arc with $k\geq 2$, and hence by Corollary~\ref{smootharc},
it has a coordinate neighborhood.

We now turn to the second conclusion. In the case $k=0$, the domain
$\Omega$ is circular, hence $\partial\Omega$ is a disjoint union of
circles. $L_\xi$ is {\em not} transverse to $\partial\Omega$ iff it
is tangent to some component circle of $\partial\Omega$. So we can
take $E$ to be the finite set of $\xi$'s such that $L_\xi$ is
tangent to $\partial\Omega$.

In the case $k\geq 1$, let $\Gamma$ be a connected component of
$\partial\Omega$. We can parameterize $\Gamma$ by a $\mathcal{C}^1$
map $\gamma=\gamma_1+i\gamma_2:S^1\rightarrow\Gamma\subset\cx$. The
line $L_\xi$ is {\em not} transverse to $\Gamma$ iff $\xi$ is a
critical value of $\gamma_2:S^1\rightarrow \rl$. By Sard's theorem,
the set $E_\Gamma$ of critical values of $\gamma_2$ is of measure 0.
Of course $E$ is closed. We let $E=\bigcup E_\Gamma$, with a union
over the finitely many components $\Gamma$ of $\partial\Omega$.
Then $E$ is nowhere dense.
\end{proof}

We make the following simple observation:
\begin{obs}\label{strip}
Let $u$ and $v$ be real valued  $\mathcal{C}^1$ functions defined on a neighborhood of $0$ in
$\rl$ such that for each $x$, we have $u(x)<0<v(x)$. Then there is an $\eta>0$
such that for $0<\theta\leq\eta$, the vertical strip
\[ S \assign  \{ (x,y)\in \rl^2\colon x\in [-\theta,\theta],
 {u}(x)\leq y \leq {v}(x) \} \]
is star shaped with respect to the origin.
\end{obs}
\begin{proof}
Clear.
\end{proof}

We will now decompose $\overline{\Omega}$ into a good pair $(K_1,K_2)$.

\begin{lemma}\label{decomposition}
 Let $F$ be as in Lemma~\ref{slicing}. There is a good pair
$(K_1,K_2)$ such that  $K_1\cup K_2= \overline{\Omega}$,
and each $F(K_j)$ has a coordinate neighborhood $\mathcal{V}_j$ in
$\cx\times\m$.
\end{lemma}
\begin{proof} Thanks to Lemma~\ref{slicing}, for each vertical line
$L$, the set $F(L\cap\overline{\Omega})$ has a coordinate
neighborhood in $\cx\times\m$. By compactness we can find finitely
many points
\[
 x_0<x_1<\cdots<x_N,
\]

such that for $j=0,\ldots, N-1$ the set $ \{F(z)\colon x_j\leq \Re
z< x_{j+1}, z\in \overline{\Omega} \} $ has a coordinate
neighborhood in $\cx\times\m$, and each component of
$\{z\in\cx\colon x_j\leq \Re z< x_{j+1}, z\in \overline{\Omega} \}$
is simply connected. We will impose the following condition on the
points $x_j$ : for each $j$, the
vertical line $\Re z=x_j$ meets $\partial \Omega$ transversely at
each point of intersection. Since the set $E$ of $\xi$'s such
that $\Re z=\xi$ is {\em not} transverse $\partial \Omega$ has been shown
above to be closed and of measure $0$, this can be easily done.

Therefore,  $\partial\Omega$ is a union of (an even number of )
graphs of $\mathcal{C}^1$ functions  in a
neighborhood of each of the vertical lines $\Re z=x_j$.
Thanks to Observation~\ref{strip} above, there is an $\eta$ such
that if $\theta\leq\eta$, each component of the intersection of
$\overline{\Omega}$ with a vertical strip of width $\theta$ about the line
$\Re z=x_j$ is star shaped.
 Define the compact subsets $K_1$ and $K_2$ of $\cx$ given by
\[
 K_1\assign \{z\in \overline{\Omega}\colon x_{2j-1}-\theta\leq\Re z
\leq x_{2j}+\theta, j=1,2,\ldots\}\] and
\[
K_2\assign \{z\in\overline{\Omega}\colon x_{2j}\leq \Re z \leq
x_{2j+1},j=1,2,\ldots \}. \] In other words, $K_2$
(resp. $K_1$) consists of the slices of $\overline{\Omega}$ over the
odd numbered intervals in the partition $x_0<x_1<\cdots<x_N$ (resp.
the slices over the even numbered ones slightly enlarged). The sets
$\overline{K_1\setminus K_2}$ and $\overline{K_2\setminus K_1}$ are
disjoint, and each $F(K_j)$ has a coordinate neighborhood in
$\cx\times\m$, which will be our $\mathcal{V}_j$.

\end{proof}

We can now use Theorem~\ref{goodpairapprox} to prove the following
approximation result:
\begin{prop} \label{pointapprox}
There is a point $p\in\partial\Omega$ with the following property:
given any $\epsilon>0$, there is a 
 neighborhood $U_\epsilon$ of $p$ in $\cx$
and a map $g\in\mathcal{A}^k\left( \overline{\Omega},\m\right)$, such that$
\tmop{dist}_{\mathcal{C}^k(\overline{\Omega},\m)}(f,g)<\epsilon$, and
$g$ extends as a holomorphic map to $U_\epsilon$.
\end{prop}
\begin{proof}
The proof is an application of Theorem~\ref{goodpairapprox}. Let as before,
$F(z)=(z,f(z))$. In the
notation of that lemma we choose the following data:
\begin{itemize}
\item the good pair $(K_1,K_2)$ will be the one in the conclusion
of Lemma~\ref{decomposition}, so that $K_1\cup K_2=\Omega$, and
each of $F(K_1)$ and $F(K_2)$ has a coordinate neighborhood in 
$\cx\times\m$.
\item Let $p\in (K_2\setminus K_1)\cap \partial\Omega$, and 
let $V$ be a closed disc around $p$ such that $V\cap K_1=\emptyset$.
 It is clear (e.g. by an easily established
$\mathcal{C}^k$ version of Mergelyan's theorem) that any function in
$\mathcal{A}^k(K_2)$ can be approximated in the $\mathcal{C}^k$
sense by entire functions, therefore a fortiori by functions in
$\mathcal{A}^k(K_2\cup V)$.
\item We will let the target manifold be $\cx\times\m$ (denoted  by
$\m$ in the statement of lemma~\ref{goodpairapprox}), and the map to
be approximated be $F$. As remarked above,
$F(K_j)$ has a coordinate neighborhood in $\cx\times\m$.
\end{itemize}
Therefore, by Lemma~\ref{goodpairapprox} we obtain a $\mathcal{C}^k$
approximation $G$ to $F$ on $K_1\cup K_2=\overline{\Omega}$, where
$G$ extends holomorphically to some neighborhood $U_\epsilon$ of $p$,
and we have $\tmop{dist}_{\mathcal{C}^k(\overline{\Omega},\m)}(F,G)<\epsilon$.
 Let $g=\pi_\m\circ G$, where
$\pi_\m:\cx\times\m\rightarrow\m$ is the projection onto the second
component.
\end{proof}

At this point, the approximation has been achieved in a neighborhood 
of one point $p$ in the boundary. We could repeat this process, thus 
obtaining a proof of Theorems~\ref{ctsapprox}$^\prime$ and 
$\ref{smoothapprox}$. This would require a strengthened version 
of Proposition~\ref{pointapprox} in which (i) we can choose the point
$p$ arbitrarily, and (ii)the diameter of the set
$\partial\Omega\cap U_\epsilon$ is independent of $p$. This is the 
route followed in \cite{thesis}.

In this paper however, we complete the proof using a technique
found in \cite{forst}, which in the setting of our problem allows
us to avoid entirely the method of  successive bumpings outlined
above. We explain this method in the following lemma.
\begin{lemma}\label{fdd}
Let $D$ be a $\mathcal{C}^1$ domain in the plane, and let $V$ be
an open set in $\cx$ such that $\partial D\cap V\not=\emptyset$. 
Then there is a one parameter family  
$\psi_t\in\mathcal{H}(\overline{D})$, such that for small $t\geq 0$,
we have $\psi_t(\overline{D})\subset D\cup V$, and as $t\rightarrow 0$,
$\psi_t\rightarrow \psi_0$ in the $\mathcal{C}^\infty$ sense, where 
$\psi_0$ is the identity map $\psi_0(z)=z$.
\end{lemma}
\begin{proof}The case in which $\partial D$ is connected being 
trivial, we assume that it has at least two components.
For a Jordan curve $C\subset \cx$ , denote by $\beta(C)$ the bounded 
component of $\cx\setminus C$. Let $\left\{ C_k\right\}_{k=1}^{M+1}$ 
be the components of $\partial D$, where $V\cap C_{M+1}\not=\emptyset$.
There is a component $C_j$ of $\partial D$ with the following 
properties :  for $k\not=j$, we have (i) $C_k\subset\beta(C_j)$ and
(ii)$\beta(C_k)\cap D =\emptyset$. Call $C_j$ the {\em outer boundary} of
$D$.

Without loss of generality we may assume that the outer boundary of $D$ is
$C_{M+1}$.  If this is not already the case, make in a neighborhood of 
$\overline{D}$ the change of coordinates $z\mapsto \frac{\rho^2}{z-z_0}$,
where $z_0\in\beta(C_{M+1})$, and $\rho>0$ is so small that
 $B_\cx(z_0,\rho)\Subset\beta(C_{M+1})$.

Let $V_0\Subset V$ be such that $V_0\cap C_{M+1}\not=\emptyset$, and
for $1\leq k\leq M$, we have $V_0\cap C_k=\emptyset$. Let $\Gamma=
\partial D\setminus V_0$. Then $\cx\setminus\Gamma$ has $M+1$ components.
Exactly one of these is unbounded, and this component contains the set 
$C_{M+1}\setminus V_0$. Let $P$ be a set of $M$  points, such that for
$j=1,\ldots,M$, each bounded component $\beta(C_j)$ of $\cx\setminus\Gamma$
contains exactly one point of $P$. Of course, $P\cap D=\emptyset$.

Let $N$ be the inward directed unit normal vector field on $\partial D$. Identifying
$T\cx$ with $\cx$, the restriction $N|_\Gamma$ is a continuous function on the set 
$\Gamma$. Since $\Gamma$ has no interior, and $\cx\setminus\Gamma$ has 
finitely many components, $N$ can be uniformly approximated on
$\Gamma$ by rational functions holomorphic on $\cx\setminus P$. We therefore obtain
a holomorphic vector field $X$ on $\overline{D}$ such that $X|_\Gamma$ is directed 
inward, i.e., towards $D$. Denote by $\psi_t$ the holomorphic flow generated by $X$.
Clearly, on any compact set, $\psi_t$ approaches the identity in all $\mathcal{C}^k$ norms
as $t\rightarrow 0$. For small $t\geq 0$, we have $\psi_t(\Gamma)\subset D$,
and by continuity, $\psi_t(\partial D\cap \overline{V_0})\subset V$. Therefore,
$\psi_t(\partial D)\subset D\cup V$, so that $\psi_t(\overline{D})\subset D\cup V$.
\end{proof}

We can now end the proof of Theorems~\ref{ctsapprox}$^\prime$ and
\ref{smoothapprox}. Suppose that $\epsilon>0$ is  given, and
we want to find a $h\in \mathcal{A}^k(\overline{\Omega},\m)$ such
that
$\tmop{dist}_{\mathcal{C}^k(\overline{\Omega},\m)}(f,h)<\epsilon$.
Using Lemma~\ref{pointapprox} above, we can construct an
approximation $g$ which extends holomorphically to a neighborhood
$U_\epsilon$ of a point $p$ on the boundary, and we have
$\tmop{dist}_{\mathcal{C}^k(\overline{\Omega},\m)}(f,g)<\frac{\epsilon}{2}$.
In lemmma \ref{fdd} let $V=U_\epsilon$, and $D=\Omega$.
Let $\psi_t$ be the family of biholomorphisms in the conclusion of
Lemma~\ref{fdd}. Then, $g_t=g\circ\psi_t$ is in
$\mathcal{H}(\overline{\Omega},\m)$ for small $t$, and as
$t\rightarrow 0$, we have $g_t\rightarrow g$ in the $\mathcal{C}^k$
sense on $\overline{\Omega}$.
 Therefore, taking $t>0$ small enough we obtain $h=g_t$ such
that
$\tmop{dist}_{\mathcal{C}^k(\overline{\Omega},\m)}(h,g)<\frac{\epsilon}{2}$.
Theorems~\ref{ctsapprox}$^\prime$ (and therefore
Theorem~\ref{ctsapprox})  and \ref{smoothapprox} are therefore
proved.

\end{document}